\documentclass[]{imsart}

\RequirePackage{amsthm,amsmath,amssymb,amsfonts}
\RequirePackage[numbers]{natbib}
\RequirePackage[colorlinks,citecolor=blue,urlcolor=blue]{hyperref}
\RequirePackage{graphicx}
\usepackage{amsfonts}
\usepackage{subfigure}

\startlocaldefs

\usepackage{mathbbol}
\DeclareSymbolFontAlphabet{\amsmathbb}{AMSb}%
\renewcommand{\P }{{\rm Pr}}
\newcommand{\E }{{\rm \mathbb{E}}}
\newcommand*\dd{\mathop{}\!\mathrm{d}}

\numberwithin{equation}{section}
\theoremstyle{plain}
\newtheorem{theorem}{Theorem}

\newtheorem{remark}[theorem]{Remark}

\endlocaldefs

\begin{document}

\begin{frontmatter}
\title{Non-parametric cure models through extreme-value tail estimation\support{M. Bladt was supported by the Carlsberg Foundation, grant CF23-1096. I. Van Keilegom gratefully acknowledges funding from the FWO and F.R.S. - FNRS (Excellence of Science programme, project ASTeRISK,
grant no. 40007517) and the FWO (senior research projects fundamental research, grant no. G047524N).}
}
\runtitle{Non-parametric cure models through extreme-value tail estimation}

\begin{aug}

\author{\fnms{Jan} \snm{Beirlant}\ead[label=e0]{jan.beirlant@kuleuven.be}}
\address{Department of Mathematics, KU Leuven, Celestijnenlaan 200B, 3001 Leuven, Belgium\\
\printead{e0}}
\author{\fnms{Martin} \snm{Bladt}\ead[label=e1]{martinbladt@math.ku.dk}}
\address{Department of Mathematical Sciences, University of Copenhagen, Denmark\\
\printead{e1}}
\author{\fnms{Ingrid} \snm{Van Keilegom}
\ead[label=e2]{ingrid.vankeilegom@kuleuven.be}}
\address{Faculty of Economics and Business, KU Leuven,
Naamsestraat 69, 3000 Leuven, Belgium\\
\printead{e2}}

\end{aug}

\begin{abstract}
In survival analysis, the estimation of the proportion of subjects who will never experience the event of interest, termed the cure rate, has received considerable attention recently. Its estimation can be a particularly difficult task when follow-up is not sufficient, that is when the censoring mechanism has a smaller support than the distribution of the target data. In the latter case, non-parametric estimators were recently proposed using extreme value methodology,  assuming that the distribution of the susceptible population is in the Fréchet or Gumbel max-domains of attraction. In this paper, we take the extreme value techniques one step further, to jointly estimate the cure rate and the extreme value index, using probability plotting methodology, and in particular using the full information contained in the top order statistics. In other words, under sufficient or insufficient follow-up, we reconstruct the immune proportion. To this end, a Peaks-over-Threshold approach is proposed under the Gumbel max-domain assumption. Next, the approach is also transferred to more specific models such as Pareto, log-normal and Weibull tail models, allowing to recognize the most important tail characteristics of the susceptible population. We establish the asymptotic behavior of our estimators under regularization. Though simulation studies, our estimators are show to rival and often outperform established models, even when purely considering cure rate estimation. Finally, we provide an application of our method to Norwegian birth registry data.
\end{abstract}

\tableofcontents
\end{frontmatter}
\section{Introduction}
In survival analysis, there is growing attention to the problem of accounting for subjects who will never experience the event of interest, known as cured or non-susceptible subjects. For this reason, the accurate estimation of the cure rate has received attention from parametric, semi-parametric and nonparametric fronts alike, cf.~\cite{PengYu2021}. However, it is well established that the asymptotic accuracy of most estimators requires sufficient follow-up, which in statistical terms means that the censoring mechanism has support which is at least as large as the support of the susceptible population. This assumption, however, is often violated in practice, and so it is difficult to determine whether an event is censored due to insufficient follow-up or due to immunity. In essence, these events can only be differentiated at high quantiles of the sample, which in turn is the subject of study of extreme value theory (EVT), and thus it is natural to approach the problem using these techniques. In this paper, we pursue this strategy.

Techniques using both survival analysis and extreme value theory have classically led to estimators which concentrate on the tail estimation of the underlying distribution, see for instance~\cite{BeirlantGuillou2001,EinmahlFilsGuillou2008,BeirlantWorms2019,BladtAlbrecherBeirlant2021}. A recent development that considers the estimation of cure rate models under max-domains of attraction conditions is provided in~\cite{EscobarKeilegom2019} and~\cite{EscobarMallerKeilegomZhao2021} for the Fr\'echet and Gumbel max-domains of attraction, respectively (see also~\cite{EscobarKeilegom2023} in a conditional setting). They propose a three-point Pickands-type estimator which corrects the nonparametric estimator of the cure rate of~\cite{MallerZhou1992}. Since the focus there is on the cure rate, in particular, no estimation of the tail of the distribution is provided. However, the general idea of their approach is intuitive: to extrapolate into the tail using EVT asymptotics, so that an additional proportion of right-censored individuals can be classified as susceptible, coming one step closer to the true cure rate. 

We adopt the general idea of using EVT to extrapolate beyond the censoring support. General reviews on EVT can be found in~\cite{Coles2001}, \cite{BGTS2004}, \cite{haan2006extreme}. We extend existing methods by using all the top order statistics, which intuitively should boost the performance of the estimator of the cure rate. This is then confirmed mathematically and in simulations. An additional benefit of our estimation method through probability plots is that we are able to recover the asymptotic tail parameters without additional effort. Thus, our model can be seen as jointly estimating the cure rate and the susceptible tail asymptotics simultaneously. In particular, our approach also provides value to the extremes field, where we are able to remove the effect of the immune population for accurate tail estimation. Ignoring the cure rate would of course suggest a much too heavy tail.
 The estimator is competitive for both sufficient and insufficient follow-up, which is a desirable property to circumvent the determination of sufficient follow-up hypotheses (see~\cite{MallerResnick2022} for sufficient follow-up testing in the Gumbel max-domain).

The remainder of the paper is structured as follows. In Section~\ref{sec:estimators} we consider the granular problem when the model falls into the the Gumbel or Fr\'echet max-domain of attraction. Subsequently in the same section we consider other tail models such as log-normal and Weibull type tail models within the Gumbel domain, next to the Pareto-type model. The asymptotic results are given in Section~\ref{sec:asymptotics}. In Section~\ref{sec:sims} we provide finite-sample behavior simulation results.  Subsequently,  the use of our method is illustrated on real data in Section~\ref{sec:real}. Finally,  Section~\ref{sec:conc} concludes. The proof of the asymptotic normality theorem is provided in the Appendix.

\section{Cure rate estimation based on extreme value methods }\label{sec:estimators}

Let the survival time of a subject be denoted by $T$ and the cure rate is $1-p$ where $p=\P(T<\infty)$ is the proportion of the population that is susceptible. Due to random right-censoring, we do not observe the survival times of all subjects. Instead we observe $Z$ and $\delta$ where $Z=\min (T,C)$,  $\delta = \mathbf{1}_{\{T \leq C \}}$  with $C$ the censoring time with distribution function (df) $G$ that is assumed to be finite, and independent of $T$. This implies that all cured individuals (with $T=\infty$) are censored, and among the susceptible population, some are censored.  The df $H(t)= \P(Z\leq t)$ satisfies
$$
1-H(t) = (1-F(t))(1-G(t)).
$$
\\
The subdistribution $F$ of $T$ is given by
\begin{equation}
F(t)= \P(T\leq t) =p F_0 (t),
\label{rate}
\end{equation}
with $F_0$ the distribution function of the survival times of the susceptible subjects. We denote by
$\hat F$ the product-limit estimator (\cite{KaplanMeier1958}) based on an independent and identically distributed (i.i.d.) sample $(Z_i, \delta_i)$, $i=1,\ldots ,n$. Let the order statistics be $Z_{1,n}\le\cdots\le Z_{n,n}$, with corresponding concominant indicators $\delta_{1,n}\le\cdots\le \delta_{n,n}$.

An important and natural estimator is given 
$p_n=\hat F(Z_{n,n})$, which was proposed as an estimator of $p$ by~\cite{MallerZhou1992}. These authors also showed that it is a consistent estimator for $p$ as $n\to \infty$ if and only $\tau_0 \leq \tau_c$ with $\tau_0$ and $\tau_c$ denoting the endpoints of the distributions of the susceptible subjects, and of the censoring mechanism, respectively. In the case of insufficient follow-up, i.e. 
$\tau_0 > \tau_c$ (and hence $\tau_c<\infty$), \cite{EscobarKeilegom2019} and \cite{EscobarMallerKeilegomZhao2021} proposed to do a simple improvement of $p_n$ considering $\hat{F}$ at $Z_{n,n}$ and at two additional points below, which results in a correction of $p_n$, though it is of course not possible to obtain consistency.  

The methods we propose are centered around the max-domain of attraction conditions. The Gumbel max-domain of attraction contains distributions for which normalized sample maxima $(T_{n,n}-b_n)/a_n$ converge in distribution to the Gumbel distribution with df $\exp(-e^{-x})$ for some appropriate sequences $a_n>0$ and $b_n$ for sample sizes $n$ tending to $\infty$. This  wide class of upper tails  contains well-known distributions such the Weibull, normal and log-normal distributions, which are of main importance in survival analysis. 

Next we also consider the Fr\'echet max-domain consisting  of power law or Pareto-type tails, hence heavier tailed than any element from the Gumbel domain.  Heavy-tailed time-to-event data are commonly encountered in reliability, information technology and finance, while being not so representative  of for instance human lifetimes. 
In our setting, Pareto-type distributions are defined through
\begin{equation}
1-F_0 (t) = t^{-1/\gamma}\ell (t),
\label{Pareto}
\end{equation}
with $\gamma >0$ termed the extreme value index, and $\ell$ a slowly varying function at infinity; i.e. for every $u>0$
$$
\frac{\ell (tu)}{\ell (t)}\to 1 \mbox{ as } t \to \infty. 
$$
An important and popular sub-class is given by the Hall-type slowly varying functions, \cite{HallWelsh1985},  assuming the  second-order condition
\begin{equation}
\ell (t) = C (1+D t^{-\beta}(1+o(1)), \; t \to  \infty,
\end{equation}
with constants $C>0$, $D$ real valued and $\beta >0$, which in practice englobes numerous distribution classes within the Fr\'echet max-domain of attraction.

\subsection{Using intermediate order statistics starting from probability plots}

The Gumbel domain definition is not quite specific and it is composed of quite different tail models ranging from Weibull up to log-normal tail decay. The log-normal tail model  is close to and hard to discriminate from the Pareto-type model. In survival analysis and reliability hence there is a need for tail model specification within the wide Gumbel class of models. In practice one can test the goodness-of-fit of specific models based on empirical probability plots. A first set of methods for estimating the cure rate is based on such approach.  We first consider the Pareto probability plot approach followed by probability plotting for Weibull and log-normal tail models. In principle, other tail behavior can be targeted, and these three approaches can be considered as archetype constructions tailored for popular distributions.

\subsubsection{Assuming Pareto-type behavior}
The Pareto-type assumption can be graphically verified by plotting
\begin{equation}
\big(\log t , \,  -\log (1- \hat{F}_0(t)) \big)
\label{Pa_PPplot}
\end{equation}
for a set of $t$ values and using some appropriate substitute $\hat{F}_0$ for $F_0$. If model \eqref{Pareto} holds, then a linear pattern emerges in \eqref{Pa_PPplot} as $t$ is taken large enough. Also note that the slope of the linear part is then an estimator of $1/\gamma$. Here we use a $Z_{n-k,n}$ value for some value of $k$ for the $t$ value and we inspect the linearity in the set of the largest $k$ observations
$$
 Z_{n-k+1,n} \le \cdots\le Z_{n,n}.
$$ 
A direct estimator of $F_0$ is not possible here but using the Kaplan-Meier estimator $\hat{F}_n$ we can use the substitute
$$
\frac{\hat{F}_n}{p} (Z_{n-j+1,n}), \; j=1,\ldots,k,
$$ 
and we perform the least-squares regression optimization to estimate the parameters $\beta$ and $p$:  
\begin{equation}\label{SS_P}
SS_p(\beta,p)= \sum_{j=1}^k\left(
-\log \frac{1-{\hat{F}_n(Z_{n-j+1,n})\over p}}{1-{\hat{F}_n(Z_{n-k,n})\over p}}
-\beta  \log {Z_{n-j+1,n} \over Z_{n-k,n}} \right)^2.
\end{equation}
The resulting estimator of $p$ is denoted by  $\hat{p}^P_k$. The slope parameter   $\beta$ leads to an estimate of $1/\gamma$.

After having obtained joint parameters for the cure rate and slope parameters, we may for a given $k$ and using only the corresponding estimate $\hat{p}_k^P$ of $p$ (and not the slope estimate), construct the plot
$$
\left(  \log Z_{n-j+1,n}, -\log (1-{\hat{F}_n(Z_{n-j+1,n})\over \hat{p}_k^P}) \right), \; j=1,\ldots,k
$$
which allows to assess the goodness-of-fit of the Pareto-type model, taking the cure rate into account. 

\subsubsection{Assessing models in the Gumbel max-domain}

The method discussed above concerning Pareto-type tails can be adapted in order to estimate $p$ and to check the fit  for more specific tail models within the Gumbel max-domain, such as Weibull and log-normal tail models. 

In case of Weibull-type tails, given by $1-F_0 (t) = e^{-\lambda t^{\tau} \ell (t)}$, a graphical check can be performed by verifying linearity at large values of $t$ in
\begin{equation}
\big(\log t,  \log\{-\log (1- {F}_0(t))\} \big).
\label{W_PPplot}
\end{equation}
As in subsection 3.1 we are then lead to least-squares optimization with respect to $\beta, p$ of
 \begin{eqnarray}
SS_w(\beta,p)= \sum_{j=1}^k\left( 
\log \frac{ \log(1-\hat F_n(Z_{n-j+1,n})/p)}{\log(1-\hat{F}_n(Z_{n-k,n})/p)}
-\beta \log {Z_{n-j+1,n} \over Z_{n-k,n}}
 \right)^2, 
 \nonumber 
\label{SS_W}
\end{eqnarray}
resulting in the  estimator  $\hat{p}^W_k$ of $p$. The slope parameter $\beta$ will provide an estimate of $\tau$.
Having jointly estimated the cure rate and slope parameters, for a given $k$ and using the corresponding estimate $\hat{p}_k^W$ of $p$, the plot
$$
\left(  \log Z_{n-j+1,n}, \log\{-\log (1-{\hat{F}_n(Z_{n-j+1,n})\over \hat{p}_k^W})\} \right), \; j=1,\ldots,k
$$
allows to assess the goodness-of-fit of the Weibull-type model. 

Similarly, a log-normal type tail can be verified on the basis of linearity at large thresholds $t$ of the plot
\begin{equation}
\big(\log t, \Phi^{-1} (F_0(t)) \big),
\label{LN_PPplot}
\end{equation}
where $\Phi^{-1}$ denotes the standard normal quantile function, which leads to minimization of 
\begin{align}
&SS_{ln}(\beta,p)\nonumber\\
&= \sum_{j=1}^k\left(
\Phi^{-1}\left({\hat{F}_n(Z_{n-j+1,n})\over p}\right)- \Phi^{-1}\left({\hat{F}_n(Z_{n-k,n})\over p}\right)
-\beta 
 \log{Z_{n-j+1,n} \over Z_{n-k,n}} \right)^2, 
\label{SS_LN}
\end{align}
with respect to $(\beta, p)$. The cure rate estimator then is denoted by $\hat{p}^{L}_n$ and the slope $\beta$ leads to estimating $1/\sigma$.
For a given $k$ and using the corresponding estimate $\hat{p}_k^L$ of $p$, the plot
$$
\left(  \log Z_{n-j+1,n}, \Phi^{-1}\left( {\hat{F}_n(Z_{n-j+1,n})\over \hat{p}_k^L}\right) \right), \; j=1,\ldots,k
$$
allows to assess the goodness-of-fit of the lognormal-type model. 

Since the parameters $p$ and $\beta$ are both  related to the regulation of the far tail of the distribution, the optimization problems can have quite flat loss surfaces for finite samples. Consequently, to prevent the estimator of $p$ to run too far away from the benchmark solution $p_n= \hat{F}_n (Z_{n,n})$, we incorporate a regularization term. We now adopt an asterisk notation, where $\ast$ is to be replaced by either of the three above models. The three penalized loss functions are then  summarized as 
\begin{align}
&SS_{*}(\beta,p) \nonumber\\
&= \sum_{j=1}^k\left(
s_*\left(1- {\hat{F}_n(Z_{n-j+1,n})\over p}\right)- s_*\left(1-{\hat{F}_n(Z_{n-k,n})\over p}\right)
-\beta_* 
 \log{Z_{n-j+1,n} \over Z_{n-k,n}} \right)^2 
\nonumber   \\
 &  \quad+\lambda (p-p_n)^2, \label{SSstar} 
\end{align}
with $\lambda >0$ and 
$$
s_* (t) = \left\{
\begin{tabular}{ll}
$\log (-\log t)$ & \mbox{ for Weibull plotting,}\\
$ \Phi^{-1}(1-t)$ & \mbox{ for log-normal plotting, } $t   \in (0,1)$, \\
$-\log t$ & \mbox{ for Pareto plotting.}
\end{tabular}
\right.
$$
The probability plots are then redefined as 
 \begin{equation}
\left(    \log Z_{n-j+1,n} ,  s_*\left( 1- {\hat{F}_n(Z_{n-j+1,n})\over \hat{p}^*_k}\right)\right),\: j=1,2,\ldots,n.
\label{pp*}
\end{equation} 

The asymptotic behavior of these cure rate estimators from these optimizations is discussed in Section~\ref{sec:asymptotics}.  The estimators $\hat{p}^*_k$ are consistent for when there is sufficient follow-up, while under insufficient follow-up they are converging to $F(\tau_c)=pF_0(\tau_c)$, unless assuming $\tau_c \to \infty$.

\subsection{Cure rate estimation using the peaks-over-threshold (POT) method}
An alternative approach to targeting specific distributions as above, is to target specific max-domains of attraction. In this section we pursue this approach for the Gumbel and Fr\'echet domains.
\subsubsection{POT estimation under the Gumbel domain}

The Gumbel domain can be characterized using  the Peaks-over-Threshold (POT) result stipulating that exceedances $T-t|T>t$ for large $t$  roughly follow the exponential law for these distributions.
More specifically
\begin{align}
\lim_{t \to \infty} \P( T-t >u\sigma(t) |t< T<\infty ) = e^{-u}, \; u>0,
\label{POT}
\end{align}
for some positive $\sigma= \sigma(t)$.
We thus focus on the exceedances 
\begin{align}
E_{k-j+1,k}=Z_{n-j+1,n}-Z_{n-k,n}, \quad j=1,\dots,k,
\end{align}
for some appropriate $k$, and  the corresponding product-limit estimator $\hat F_k$ of these exceedances.
Minimizing the penalized square loss function
\begin{align}
SS_E(\sigma,\pi)= \sum_{j=1}^k\left( E_{k-j+1,k}  +\sigma \log(1-{\hat F_k(E_{k-j+1,k})\over \pi})\right)^2 +\lambda (p-p_n)^2,
\label{SS_E}
\end{align}
we then obtain $(\hat{\sigma}_k,\hat{\pi}_k)$ with $E_{k-j+1,k}$ denoting the $j$-th largest exceedance, leads to the estimator $\hat{p}^G_n$ of $p$ defined from
 \begin{align}
1-\hat{\pi}_k = \frac{1-\hat{p}^G_n}{\hat{p}(k)}
\end{align}
where $\hat{p}(k)=1-\hat{F}_n(Z_{n-k,n})$.   Note that in $\hat{p}(k)$  we are using the original product-limit estimator $\hat{F}_n$ of the $Z$ observations, rather than of the exceedances.

\subsubsection{POT estimation for Pareto-type distributions}

In this case  a  POT approach similar to the solution in the preceding subsection can be used based however on relative exceedances. Indeed, an exceedance interpretation of Pareto-type distributions is that the conditional distribution of $T/t$ given that $T>t$ for a large threshold $t$ can be approximated by the simple Pareto model: 
\begin{equation}
\bar{F}_{0,t}(y):= {\bar{F}_0(ty) \over \bar{F}_0(t)}= \P\left( T /t > y |t< T<\infty\right) \to_{t \to \infty}  
y^{-1/\gamma}
, \; y>1.
\label{POT_P}
\end{equation}
Then, using the random threshold $Z_{n-k,n}$, the log-exceedances $$E^+_{k-j+1,k}= \log (Z_{n-j+1,n}/Z_{n-k,n})= \log Z_{n-j+1,n} - \log Z_{n-k,n}$$ approximately follow the exponential distribution with mean $\gamma$, and  we minimize a similar loss function as in \eqref{SS_E} with respect to $(\gamma,\pi)$:
\begin{align}
SS_{E^+}(\gamma,\pi)= \sum_{j=1}^k\left( E^+_{k-j+1,k}  +\gamma \log(1- \frac{\hat{F}_k(E^+_{k-j+1,k})}{\pi})\right)^2 +\lambda (p-p_n)^2, 
\label{SS_Pa_POT}
\end{align}
with $\hat{F}_k$ denoting the product-limit estimator of the log-exceedances $E^+$. The resulting estimator of $p$, denoted by  $\hat{p}^F_n$, now follows from \begin{align}
1-\hat{\pi}_k = \frac{1-\hat{p}^F_n}{\hat{p}(k)}.
\end{align}

\section{Asymptotic theory}\label{sec:asymptotics}

We next consider asymptotic results under insufficient follow-up for the estimators $\hat{p}_k^*$  with $*$ referring to the Pareto (P), Weibull (W) or log-normal (L) case as proposed in Section~\ref{sec:estimators}. To this end
 we assume that the censoring distribution $G$ belongs to the Weibull max-domain of attraction with extreme value index $\gamma_c <0$:
$$
1-G(x) = (\tau_c-x)^{-1/\gamma_c},\; x < \tau_c,
$$ 
while $F_0$ is assumed to be one of the three cases 
$$
-\log \bar{F}_{0,*} (x) = 
\begin{cases}
A_w x^{\tau}\left(1+B_w x^{-\beta_w}(1+o(1))\right) & \mbox{for Weibull plotting,} \\
\log^2 x/(2\sigma^2) +A_l +B_l (\log x)^{-\beta_l}  & \mbox{for log-normal plotting, } \\
 \log (x^{\gamma^{-1}}/A_p) -B_p  x^{-\beta_p}(1+o(1))  & \mbox{for Pareto plotting,}
\end{cases}
$$
as $x \to \infty$, with $A_*, B_*$ denoting real constants, $\beta_*>0$.
It then follows that 
\begin{equation}
s_* (\bar{F}_{0,*}(x)) = C_* + \beta_* \log x
+ D_* x^{-\nu_*}(1+o(1))
\label{sF0}
\end{equation}
with $C_*, D_*$ denoting real constants, $\nu_*>0$ and  
$$
\beta_* = \left\{
\begin{tabular}{ll}
$\tau$ & \mbox{ for Weibull plotting,}\\
$1/\sigma $ & \mbox{ for log-normal plotting, }  \\
$1/\gamma$ & \mbox{ for Pareto plotting.}
\end{tabular}
\right.
$$
the slope parameter appearing in \eqref{SSstar}.

Now the distribution function of the censored data is given by 
\begin{align*}
1-H(x) &= (\tau_c-x)^{-1/\gamma_c} \times \left( 1-p +p \bar{F}_{0,*}(x)\right)\\
&=
(\tau_c-x)^{-1/\gamma_c} \times
\left[ 1-p +p \bar{F}_{0,*}(\tau_c) +pf_{0,*}(\tau_c)(\tau_c -x)(1+o(1))\right] \\
&= (\tau_c-x)^{-1/\gamma_c} \times 
\left[ 1-p_0(\tau_c) +pf_{0,*}(\tau_c)(\tau_c -x)(1+o(1))\right], \; x  \to \tau_c,
\end{align*}
with $p_0(\tau_c) = pF_{0,*}(\tau_c)$ and $f_{0,*}$ denoting the density of $F_{0,*}$.
Concerning the quantile function $Q_H$ of  $H$ one finds that
\begin{align}
 U_H (y) & :=  Q_H (1-y^{-1}) \nonumber \\
&= \tau_c - y^{\gamma_c}(1-p_0(\tau_c))^{\gamma_c}\left(1+ {1 \over 4}(1-p_0(\tau_c))^{\gamma_c-1}pf_{0,*}(\tau_c)\gamma_c y^{\gamma_c}(1+o(1)) \right)
 \label{LR}
\end{align}
as $y \to \infty$.

Next for every $y$ we have that
\begin{align}
1-{\hat{F}_n(y) \over \hat{p}^*} &=  
1- {p_0(\tau_c) \over \hat{p}^*} {\hat{F}_n(y) \over p_0(\tau_c)} \nonumber \\
&=
\bar{F}_{0,*}(y) -\left( {p_0(\tau_c) \over \hat{p}^*}-1 \right) {\hat{F}_n(y) \over p_0(\tau_c)} 
-{1 \over p_0(\tau_c)}\left( \hat{F}_n(y) - F(y)\right)
\label{decomp}
\\
& 
-{F_{0,*}(y) \over p_0(\tau_c)}\left(p-p_0(\tau_c)\right). \nonumber
\end{align}
The last term $p-p_0(\tau_c) = p\bar{F}_{0,*}(\tau_c)$ in  \eqref{decomp}
leads to  bias terms which are smaller for lighter tailed     distributions $\bar{F}_{0,*}$  such as for Weibull distributed data compared to Pareto data.   

Concerning the term in $\hat{F}_n(y) - F(y)$,
Theorem 3.14 in~\cite{MallerZhou1996} states that 
the empirical Kaplan-Meier process
$$
Y_n(t) = \sqrt{n} \big(\hat{F} (t) - F(t) \big), \; t \in [0,\tau_c]
$$
converges weakly to $\mathbf{D}(t)=(1-F(t))\mathbf{Z}(t)$ on $t\in [0,\tau_c]$ as $\tau_c < \tau_0$,  where $\mathbf{Z}$ is a Gaussian process with independent increments, mean 0, and variance process
$$
v(t)= \int_0^t \frac{\dd F(s)}{\{1-F(s)\}\{1-F(s^-)\}\{1-G(s))\}}.
$$

In order to state our main asymptotic result we introduce some further notation:
\begin{align*}
T_{k,n}&= {1 \over \sqrt{n}k }(1-p_0(\tau_c)) \sum_{j=1}^k
\left(1- \left({j \over k+1} \right)^{-\gamma_c} \right)\\
&\quad\times[\mathbf{Z} (U_H({n+1 \over j}))- \mathbf{Z}(U_H({n+1 \over k+1}))],\\
A(\tau_c) &= -1 + F_{0,*}(\tau_c)(s''_*/s'_*)(\bar{F}_{0,*}(\tau_c)), \\
M(\tau_c) &= \tau_c (1-p_0(\tau_c))^{-\gamma_c}{(1-\gamma_c)(1-2\gamma_c)\over 2\gamma_c^2},\\
B_v&=p\gamma_c (1-p_0(\tau_c))^{\gamma_c-1} f_{0,*}(\tau_c),\\
h_{1+\gamma_{c}}(t) &= \int_1^t u^{\gamma_c}\dd u,
\end{align*}
\begin{align*}
{\bf I}& =  
\begin{pmatrix}
-1& & {A(\tau_c) \over p F^2_{0,*}(\tau_c) }(\beta_*-D_*\nu_*\tau_c^{-\nu_*}) \\ \\
-A(\tau_c) & &
{A^2(\tau_c) \over p F^2_{0,*}(\tau_c) }(\beta_*-D_*\nu_*\tau_c^{-\nu_*})
- C_{\lambda} {M(\tau_c)\tau_c \over p (1-p_0(\tau_c))^{\gamma_c}(\beta_*-D_*\nu_*\tau_c^{-\nu_*})}
\end{pmatrix}.
\end{align*}

\begin{theorem}[Asymptotic representation]\label{th:asymptotic}
{ 
Assume that $\lambda=\lambda_{k,n}= C_{\lambda}\left(  {k \over n}\right)^{-2\gamma_c}$ for some $C_{\lambda} >0$, and 
$k \,n^{\gamma_c/(1-\gamma_c)} \to \infty$, then we have the following asymptotic distributional identity
\begin{align*}
& 
\begin{pmatrix}
\hat{\beta}_*-(\beta_*-D_*\nu_*\tau_c^{-\nu_*}) \\ \hat{p}^*_k- p_0(\tau_c)
\end{pmatrix} \stackrel{d}{=} M(\tau_c)(1+o_p(1)) \\
&\times\, {\bf I}^{-1}
\begin{pmatrix}
{s'_* (\bar{F}_{0,*}(\tau_c))\over p_0(\tau_c)} (n/k)^{-\gamma_c} T_{k,n} \\
{A(\tau_c)(n/k)^{-\gamma_c} T_{k,n}s'_* (\bar{F}_{0,*}(\tau_c))\over p_0(\tau_c)} - C_{\lambda}n^{-1/2}\mathbf{Z}(U_H(n)){M(\tau_c)\tau_c (1-p_0(\tau_c))^{1-\gamma_c} \over p (\beta_*-D_*\nu_*\tau_c^{-\nu_*})}
\end{pmatrix} \\
& + {\bar{F}_{0,*}(\tau_c) \over F_{0,*}(\tau_c) }A(\tau_c) (\beta_* -D_*\nu_*\tau_c^{-\nu_*}){\bf I}^{-1}
\begin{pmatrix} 1 \\ A(\tau_c) \end{pmatrix}(1+o(1)) \\
&= 
\left(\begin{pmatrix}
{-(n/k)^{-\gamma_c} T_{k,n} s'_* (\bar{F}_{0,*}(\tau_c))\over p_0(\tau_c)} \left[ M(\tau_c)+ {A(\tau_c)(\beta_* -D_*\nu_*\tau_c^{-\nu_*})^2 (1-p_0(\tau_c))^{\gamma_c} \over \tau_c C_{\lambda}F^2_{0,*}(\tau_c)}\right]\\
n^{-1/2}\mathbf{Z}(U_H(n))\; (1-p_0(\tau_c))M(\tau_c) 
\end{pmatrix}\right. \\
\\
& + \left.\begin{pmatrix} -{\bar{F}_{0,*}(\tau_c) \over F_{0,*}(\tau_c) }A(\tau_c) (\beta_* -D_*\nu_*\tau_c^{-\nu_*})
 \\ 0 \end{pmatrix}\right)(1+o_p(1)).
\end{align*}
Furthermore, $(n/k)^{-\gamma_c} T_{k,n} \stackrel{d}{=} (n/k)^{-\gamma_c/2}k^{-1/2}(1+o(1)) N(0,B_v (1-p_0(\tau_c))^{2} 
\sigma^2_k)$ 
with 
$$\sigma^2_k=
{1 \over k^2}
\sum_{j_1=1}^k \sum_{j_2=1}^k 
\left(1- \left( {j_1 \over k+1}\right)^{-\gamma_c} \right) 
\left(1- \left( {j_2 \over k+1}\right)^{-\gamma_c} \right)
 h_{1+\gamma_c}\left( {k+1 \over j_1 \vee j_2} \right),
 $$
and $n^{-1/2}\mathbf{Z}(U_H(n))\stackrel{d}{=} \sqrt{h_{1+\gamma_c}(n)/n}(1+o(1))\, N(0,1)$, are asymptotically uncorrelated.

}
\end{theorem}

\begin{remark}\rm
We provide several comments which follow from the above theorem. 
\begin{enumerate}
\item It is immediately seen that when $\lambda_{k,n} (k/n)^{-2\gamma_c}\to 0$, the matrix ${\bf I}$ is not invertible. This explains the need for regularization.
\item The condition $k \,n^{\gamma_c/(1-\gamma_c)} \to \infty$ guarantees that $(n/k)^{-\gamma_c}T_{k,n} \stackrel{\P}{\to} 0$ and hence that $\hat{\beta}_* -(\beta_*-D_*\nu_*\tau_c^{-\nu_*}) \stackrel{\P}{\to} 0.$
\item When $k=M_k  n\{h_{1+\gamma_c}(n)\}^{1/(\gamma_c-1)}$ with $M_k>0$ bounded, we have that the two random terms $(n/k)^{-\gamma_c}T_{k,n}$ and 
$n^{-1/2}\mathbf{Z}(U_H(n))$ are of the same asymptotic order.
\item Note that the bias of $\hat{p}^*_k$ has two origins. First the bias term $p_0(\tau_c)-p$ can only decrease with larger values of $\tau_c$. Next, with extreme value methods bias arises  from the second order term assumptions (containing the constants $B_*$ and $\beta_*$) in the expressions for $\bar {F}_{0,*}$ at the start of this section. Note however that the above result states that this bias term asymptotically disappears in the estimation of $\hat{p}^*_k$.
\end{enumerate}
\end{remark}

\section{Finite-sample behavior}\label{sec:sims}

In this section we investigate through simulation the finite-sample performance of the estimators defined for Pareto and Gumbel tails given by the optimization of the loss functions of equations  \eqref{SS_E}, \eqref{SS_Pa_POT} and the three models comprised in \eqref{SSstar}, where appropriate. The regularization parameter $\lambda$ in each case is taken as $k/n$. 

We consider the following \emph{scenarios}:
\begin{enumerate}
\item The exponential distribution, with df $F(x)=1-\exp(-\nu_1 x)$, $x\ge1$, $\nu_1=1$, and sufficient follow-up, that is with $G(x)=1-\exp(-\nu (x-1))$, $x\ge1$, $\nu=1/20$. The sample fraction is taken as $k=n-1$.
\item The exponential distribution, with df $F$ as in the previous case, but with insufficient follow-up with $G$ the df of a uniform random variable on $(0,3)$. The sample fraction is taken as $k=n-1$.
\item The standard lognormal distribution and insufficient follow-up with $G$ the df of a uniform random variable on $(0,6)$. The sample fraction is taken as $k=\lfloor n/5 \rfloor$.
\item The standard lognormal distribution and more insufficient follow-up with $G$ the df of a uniform random variable on $(0,2)$. The sample fraction is taken as $k=\lfloor n/5 \rfloor$.
\item The Weibull distribution, with df $F(x)=1-\exp(-x^a)$, $x\ge0$, $a=0.5$, and insufficient follow-up with $G$ the df of a uniform random variable on $(0,6)$. The sample fraction is taken as $k=\lfloor n/5 \rfloor$.
\item The Weibull distribution with df $F$ as in the previous case and more insufficient follow-up with $G$ the df of a uniform random variable on $(0,2)$. The sample fraction is taken as $k=\lfloor n/5 \rfloor$.
\item The Pareto distribution, with df $F(x)=1-x^{-1/\gamma}$, $x\ge1$, $\gamma=0.5$, and sufficient (but rather lighter-tailed) follow-up, that is with $G(x)=1-\exp(-\nu (x-1))$, $x\ge1$, $\nu=1/20$. The sample fraction is taken as $k=n-1$.
\item The Pareto distribution, with df $F$ as in the previous case, but with insufficient follow-up with $G$ the df of a uniform random variable on $(1,5)$. The sample fraction is taken as $k=n-1$.
\item The Burr distribution, with df $F(x)=1-(1+x^c)^{-d}$, $x\ge0$, $c=d=3/2$ and insufficient follow-up with $G$ the df of a uniform random variable on $(0,4)$. The sample fraction is taken as $k=\lfloor n/5 \rfloor$.
\item The Burr distribution, with df $F$ as the previous case and more insufficient follow-up with $G$ the df of a uniform random variable on $(0,2)$. The sample fraction is taken as $k=\lfloor n/5 \rfloor$.
\end{enumerate}


For all the above cases, the sample size considered is $n=5000$, and we simulate $500$ samples each time. We consider two possibilities: $p=0.9$ and $p=0.95$. 

The non-parametric estimators that we compare against our own are the following \emph{benchmarks}:
\begin{enumerate}
\item[i)] $p_n=\hat F(Z_{n,n})$ from~\cite{MallerZhou1992}.
\item[ii)] $p_y$ from~\cite{EscobarKeilegom2019}, implementing $y$ through the bootstrap procedure described in their Section 4, except that we set $\mathcal{H}=\{0.02,0.04,\dots,0.98\}$\footnote{The smaller set $\mathcal{H}=\{0.6,0.62,\dots,0.98\}$ proposed in~\cite{EscobarKeilegom2019} in general can give non-feasible correction factors.}, for distributions in the Fr\'echet max-domain of attraction.
\item[ii)] $p_G(n,\varepsilon)$ from~\cite{EscobarMallerKeilegomZhao2021}, implementing $\varepsilon$ through the least squares procedure described in their section 4, for distributions in the Gumbel max-domain of attraction.
\end{enumerate}

\begin{figure}[!htbp]
\centering
\includegraphics[width=0.24\textwidth, trim= 0in 0.25in 0.1in 0in,clip]{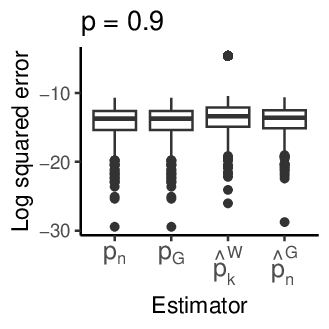}
\includegraphics[width=0.24\textwidth, trim= 0in 0.25in 0.1in 0in,clip]{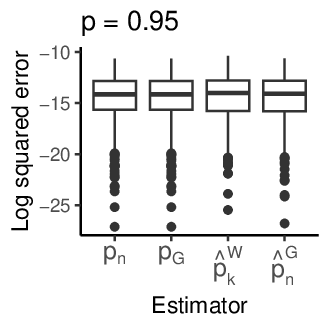}
\includegraphics[width=0.24\textwidth, trim= 0in 0.25in 0.1in 0in,clip]{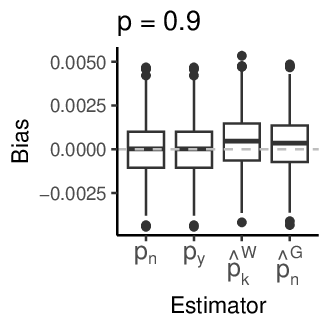}
\includegraphics[width=0.24\textwidth, trim= 0in 0.25in 0.1in 0in,clip]{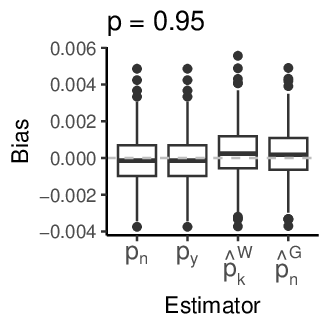}

\includegraphics[width=0.24\textwidth, trim= 0in 0.25in 0.1in 0in,clip]{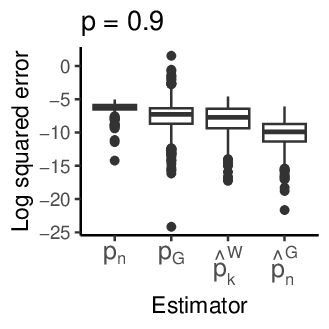}
\includegraphics[width=0.24\textwidth, trim= 0in 0.25in 0.1in 0in,clip]{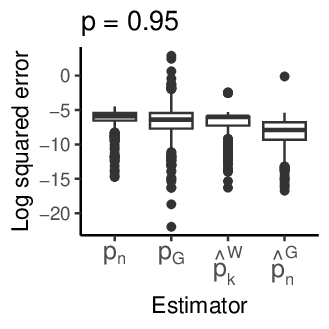}
\includegraphics[width=0.24\textwidth, trim= 0in 0.25in 0.1in 0in,clip]{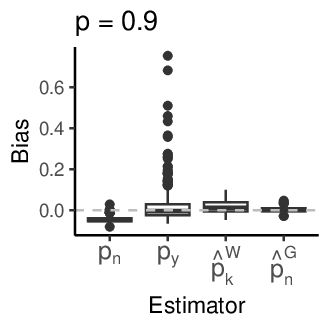}
\includegraphics[width=0.24\textwidth, trim= 0in 0.25in 0.1in 0in,clip]{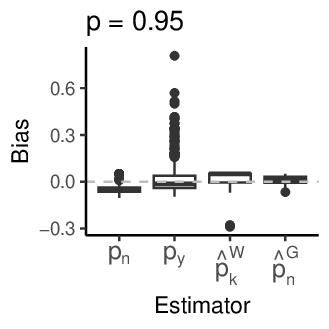}

\includegraphics[width=0.24\textwidth, trim= 0in 0.25in 0.1in 0in,clip]{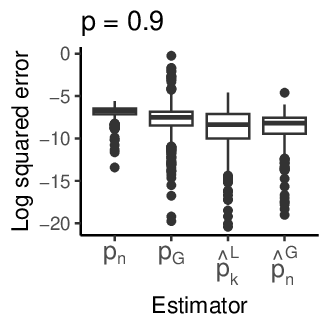}
\includegraphics[width=0.24\textwidth, trim= 0in 0.25in 0.1in 0in,clip]{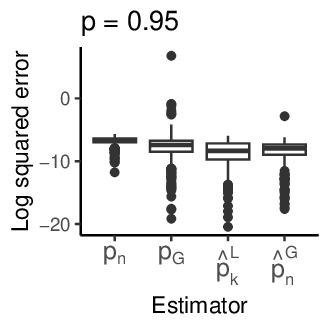}
\includegraphics[width=0.24\textwidth, trim= 0in 0.25in 0.1in 0in,clip]{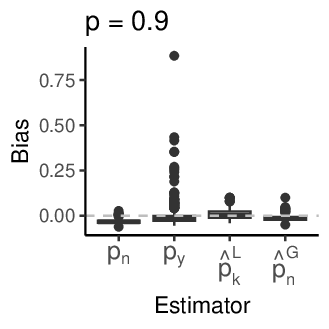}
\includegraphics[width=0.24\textwidth, trim= 0in 0.25in 0.1in 0in,clip]{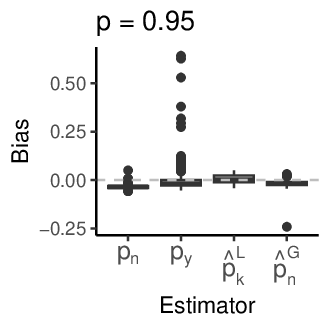}

\includegraphics[width=0.24\textwidth, trim= 0in 0.25in 0.1in 0in,clip]{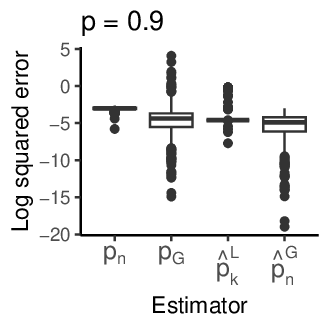}
\includegraphics[width=0.24\textwidth, trim= 0in 0.25in 0.1in 0in,clip]{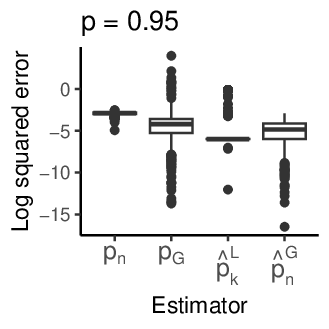}
\includegraphics[width=0.24\textwidth, trim= 0in 0.25in 0.1in 0in,clip]{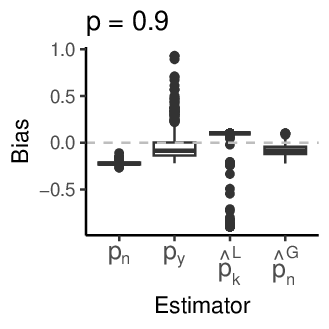}
\includegraphics[width=0.24\textwidth, trim= 0in 0.25in 0.1in 0in,clip]{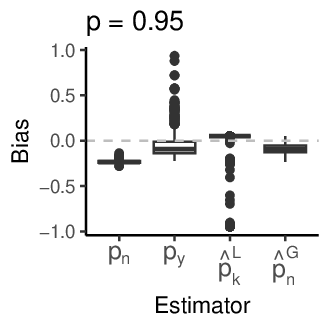}

\includegraphics[width=0.24\textwidth, trim= 0in 0.25in 0.1in 0in,clip]{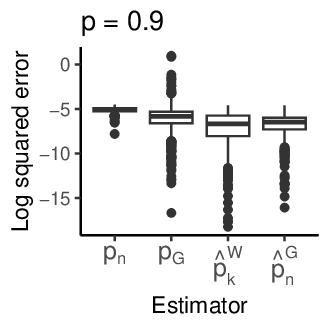}
\includegraphics[width=0.24\textwidth, trim= 0in 0.25in 0.1in 0in,clip]{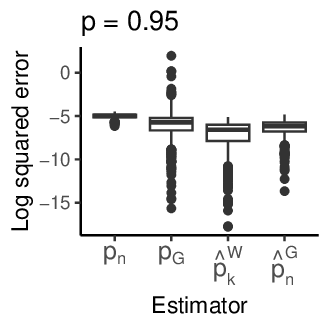}
\includegraphics[width=0.24\textwidth, trim= 0in 0.25in 0.1in 0in,clip]{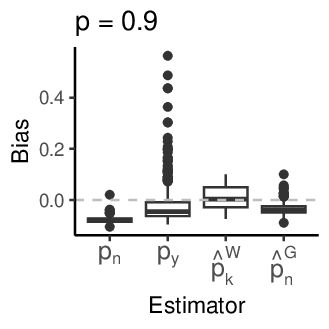}
\includegraphics[width=0.24\textwidth, trim= 0in 0.25in 0.1in 0in,clip]{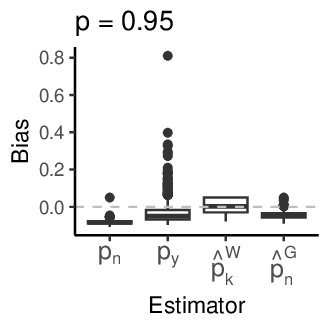}

\caption{Simulation results for scenarios 1--5 (from top to bottom).}
\label{fig:scens1to5}
\end{figure}

\begin{figure}[!htbp]
\centering

\includegraphics[width=0.24\textwidth, trim= 0in 0.25in 0.1in 0in,clip]{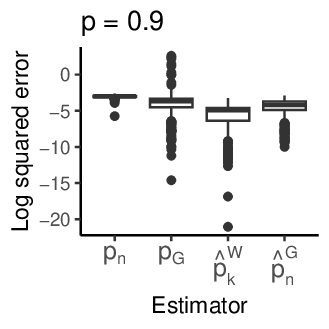}
\includegraphics[width=0.24\textwidth, trim= 0in 0.25in 0.1in 0in,clip]{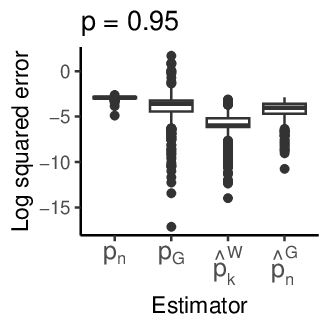}
\includegraphics[width=0.24\textwidth, trim= 0in 0.25in 0.1in 0in,clip]{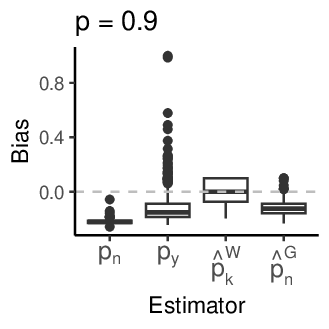}
\includegraphics[width=0.24\textwidth, trim= 0in 0.25in 0.1in 0in,clip]{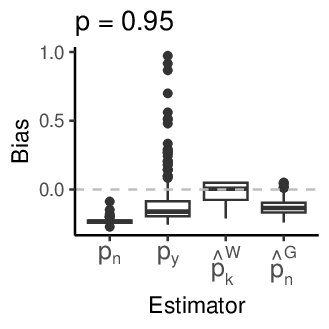}

\includegraphics[width=0.24\textwidth, trim= 0in 0.25in 0.1in 0in,clip]{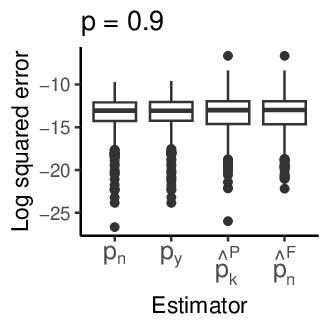}
\includegraphics[width=0.24\textwidth, trim= 0in 0.25in 0.1in 0in,clip]{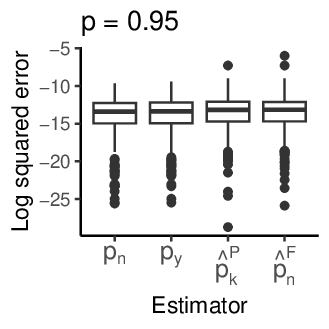}
\includegraphics[width=0.24\textwidth, trim= 0in 0.25in 0.1in 0in,clip]{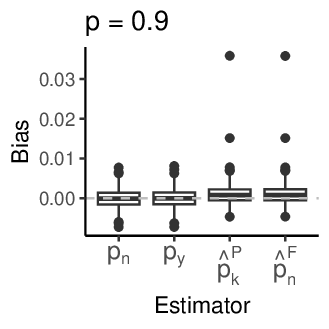}
\includegraphics[width=0.24\textwidth, trim= 0in 0.25in 0.1in 0in,clip]{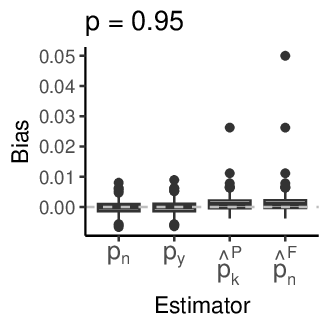}

\includegraphics[width=0.24\textwidth, trim= 0in 0.25in 0.1in 0in,clip]{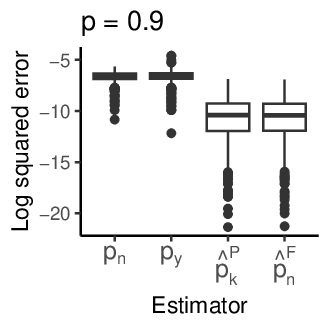}
\includegraphics[width=0.24\textwidth, trim= 0in 0.25in 0.1in 0in,clip]{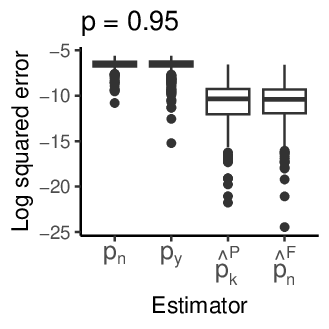}
\includegraphics[width=0.24\textwidth, trim= 0in 0.25in 0.1in 0in,clip]{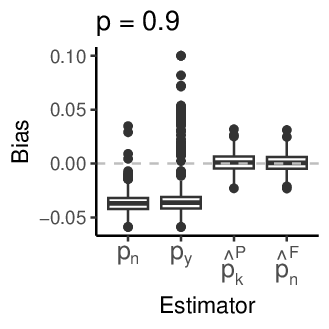}
\includegraphics[width=0.24\textwidth, trim= 0in 0.25in 0.1in 0in,clip]{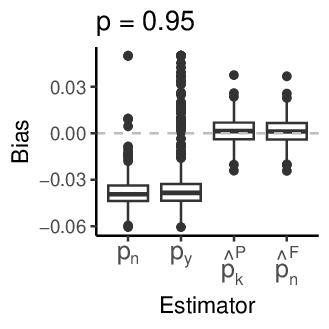}

\includegraphics[width=0.24\textwidth, trim= 0in 0.25in 0.1in 0in,clip]{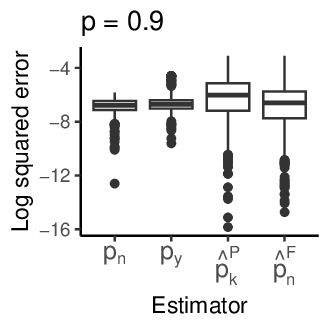}
\includegraphics[width=0.24\textwidth, trim= 0in 0.25in 0.1in 0in,clip]{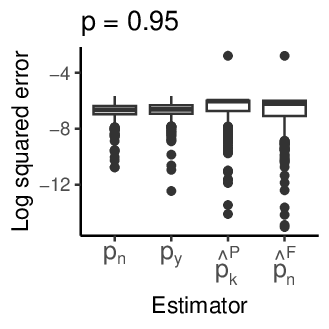}
\includegraphics[width=0.24\textwidth, trim= 0in 0.25in 0.1in 0in,clip]{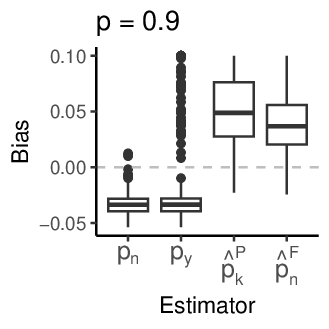}
\includegraphics[width=0.24\textwidth, trim= 0in 0.25in 0.1in 0in,clip]{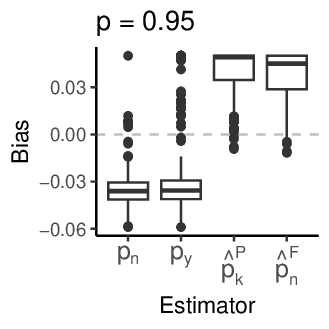}

\includegraphics[width=0.24\textwidth, trim= 0in 0.25in 0.1in 0in,clip]{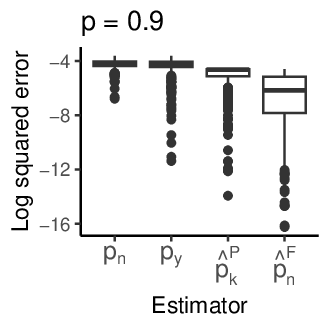}
\includegraphics[width=0.24\textwidth, trim= 0in 0.25in 0.1in 0in,clip]{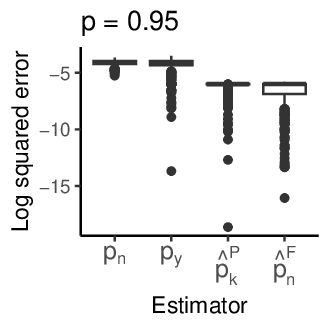}
\includegraphics[width=0.24\textwidth, trim= 0in 0.25in 0.1in 0in,clip]{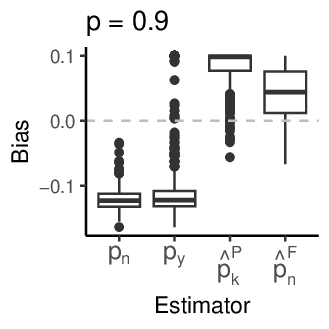}
\includegraphics[width=0.24\textwidth, trim= 0in 0.25in 0.1in 0in,clip]{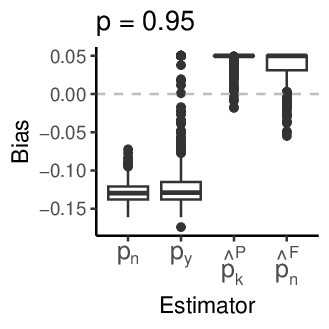}
\caption{Simulation results for scenarios 6--10 (from top to bottom).}
\label{fig:scens6to10}
\end{figure}

The distribution of the squared error losses and of the bias from the $500$ simulations is provided in Figures \ref{fig:scens1to5} and \ref{fig:scens6to10}. We observe that the extreme-value correction is not particularly favorable or unfavorable when there is a sufficient follow-up, as expected. Nonetheless, the proposed estimators $\hat p$ still may slightly outperform the other estimators in that case. In the other cases, except for  Scenario 9, the proposed  approximations provide better estimates than the state-of-the-art estimators for insufficient follow-up, in terms of squared errors. We believe this is because of the more efficient extreme-value approximation based on all $k$ upper order statistics. We also observe that the bias terms have mostly comparable behavior across estimators, with our proposed methods tending to have an upward rather than downward bias in difficult cases.
%

\section{Real data analysis}\label{sec:real}

We analyze data from the Norwegian medical birth registry. This data set is comprised of $n=53,558$ observations, which can be used to study the time between 1st and 2nd birth, with the obvious cure in this case corresponding to mothers having ultimately only one child. Whether a mother is ``cured" from having a second child, or merely right-censored (and could have a second child in the future) is the delicate feature we are trying to disentangle from the estimation procedure.

\begin{figure}[!htbp]
\centering
\includegraphics[width=\textwidth, trim= 0in 0in 0in 0in,clip]{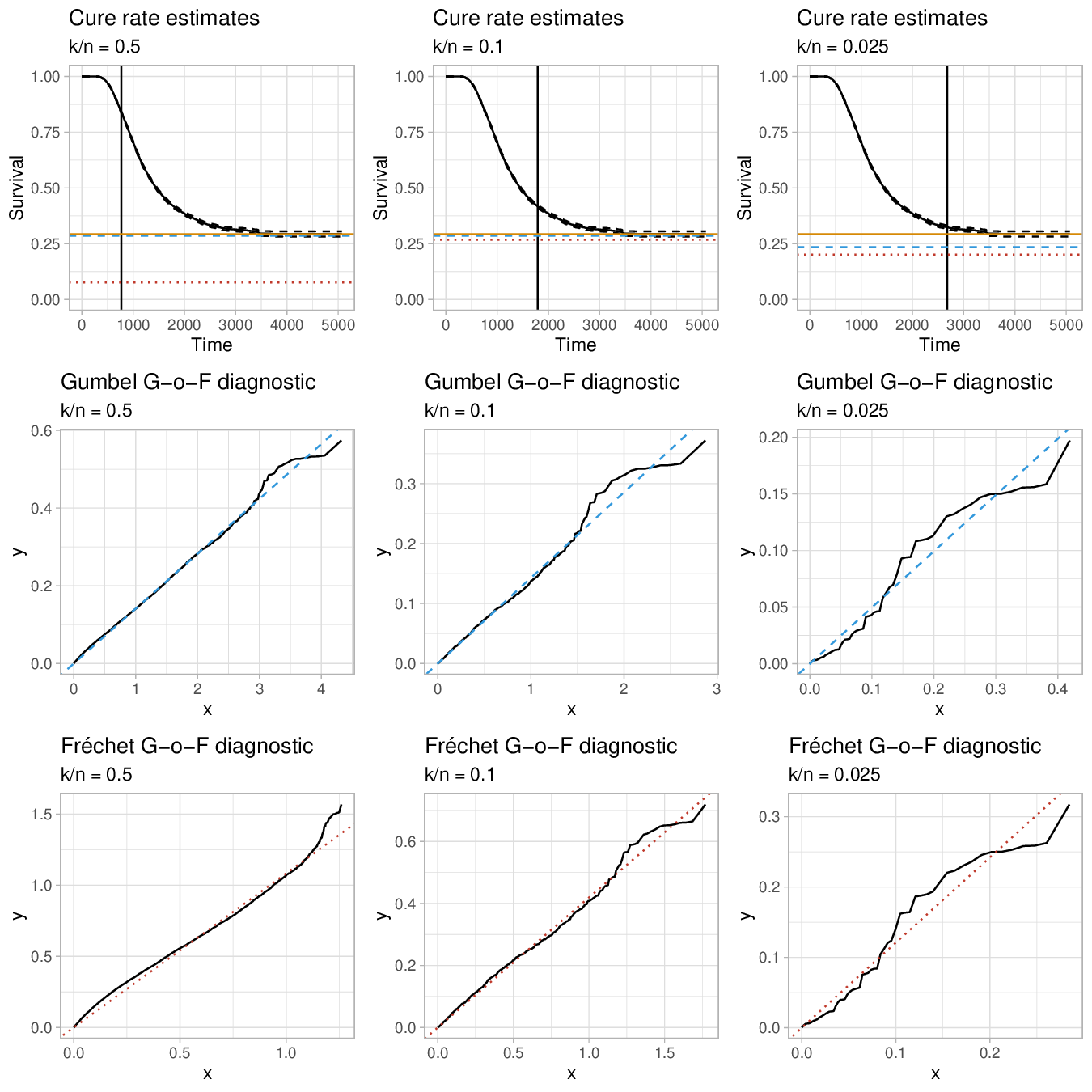}
\caption{Norwegian second borns. The $p_n$, $\hat{p}_k^G$ and $\hat{p}_k^F$ estimates, next to the Gumbel and Fr\'echet goodness-of-fit plots.  }
\label{fig:norwegian}
\end{figure}

\begin{figure}[!htbp]
\centering
\includegraphics[width=\textwidth, trim= 0in 0in 0in 0in,clip]{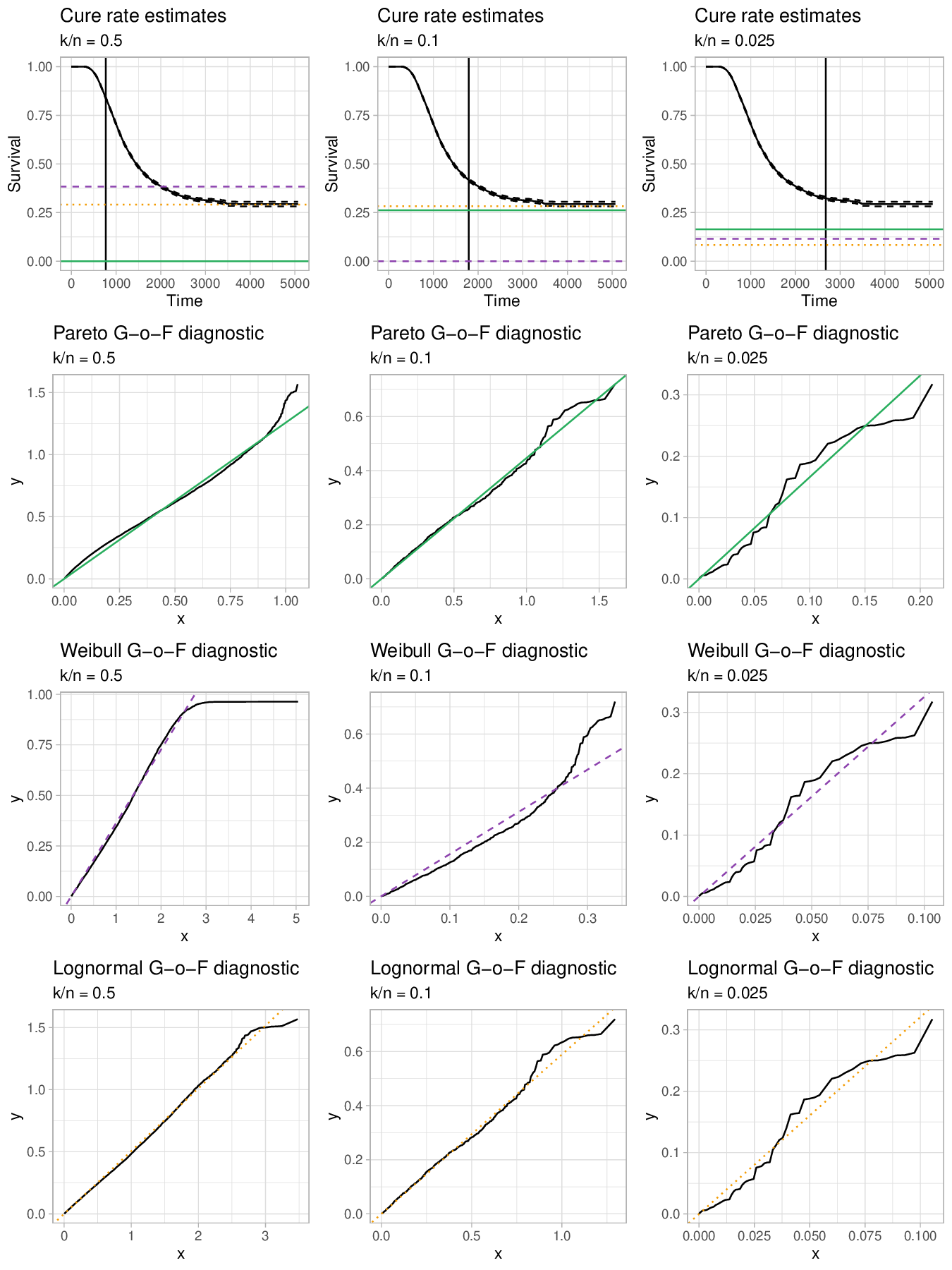}
\caption{Norwegian second borns. The $p_n$, $\hat{p}_k^P$, $\hat{p}_k^W$ and $\hat{p}_k^L$ estimates, next to the Pareto, Weibull and lognormal goodness-of-fit plots.}
\label{fig:norwegian2}
\end{figure}

In Figure \ref{fig:norwegian} we propose the Kaplan-Meier survival function, the cure rate estimates $\hat{p}_k^G$ and $\hat{p}_k^F$ jointly with the corresponding goodness-of-fit plots, each at $k/n=0.5, 0.1$ and $0.025$. The Gumbel plot shows the better fit at the large $k$ while for smaller $k$ no clear favorable model appears. In Figure \ref{fig:norwegian2} in a similar way we present the $\hat{p}_k^P$, $\hat{p}_k^W$ and $\hat{p}_k^L$ estimates together with the corresponding goodness-of-fit plots. Here at larger $k$ the lognormal goodness-of-fit plot appears to be most linear, while again for smaller $k$ the differences are less prominent. The corresponding estimates are $\hat{p}_{n/2}^G=0.714$ and $\hat{p}_{n/2}^L=0.709$, while $p_n=0.707$. 

\begin{figure}[!htbp]
\centering
\includegraphics[width=\textwidth, trim= 0in 0in 0in 0in,clip]{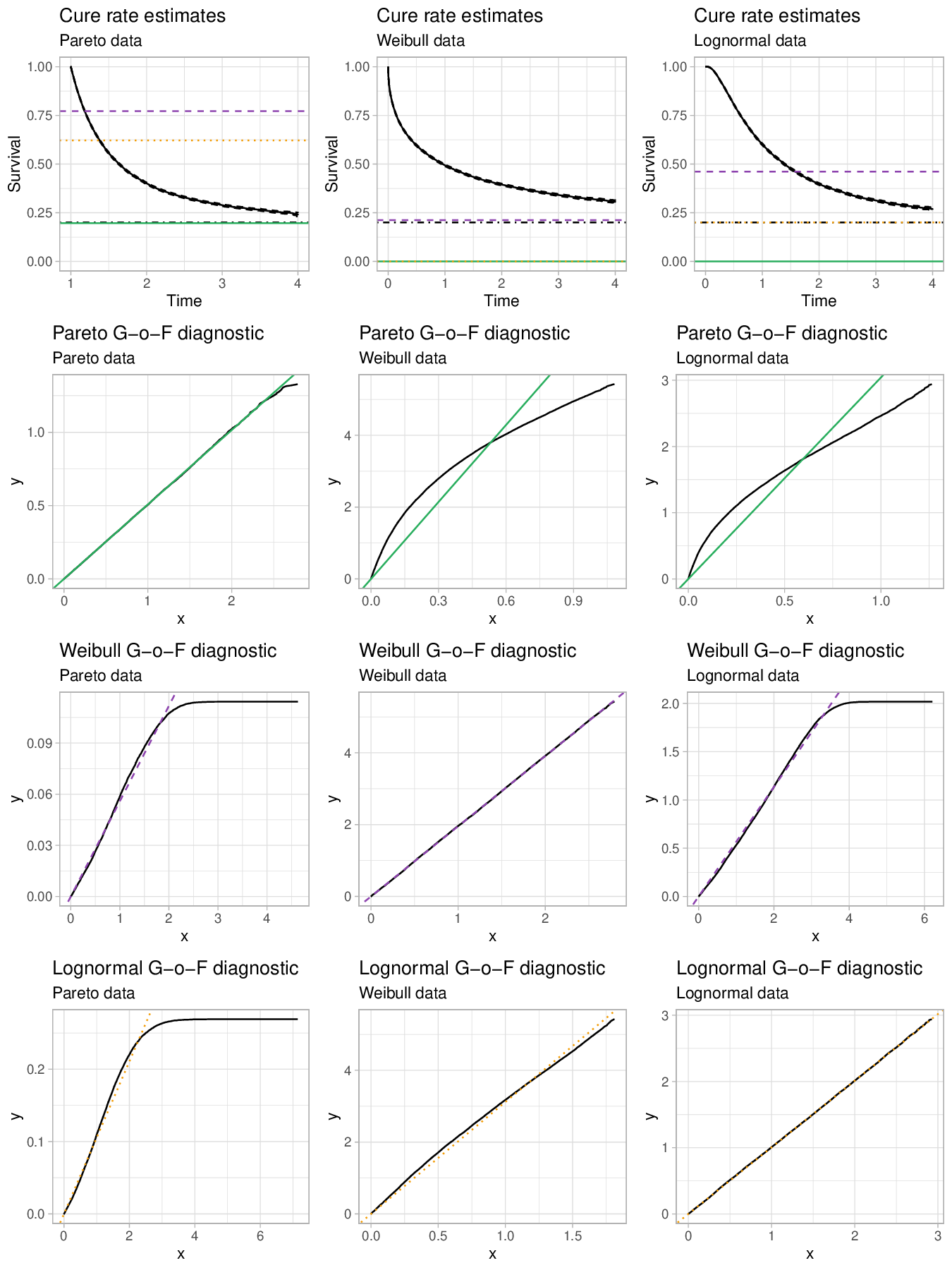}
\caption{Simulated data following the  Pareto, Weibull and lognormal models, best fitting  to the  Norwegian data. Here $k/n=0.9$}
\label{fig:sims}
\end{figure}

\begin{figure}[!htbp]
\centering
\includegraphics[width=\textwidth, trim= 0in 0in 0in 0in,clip]{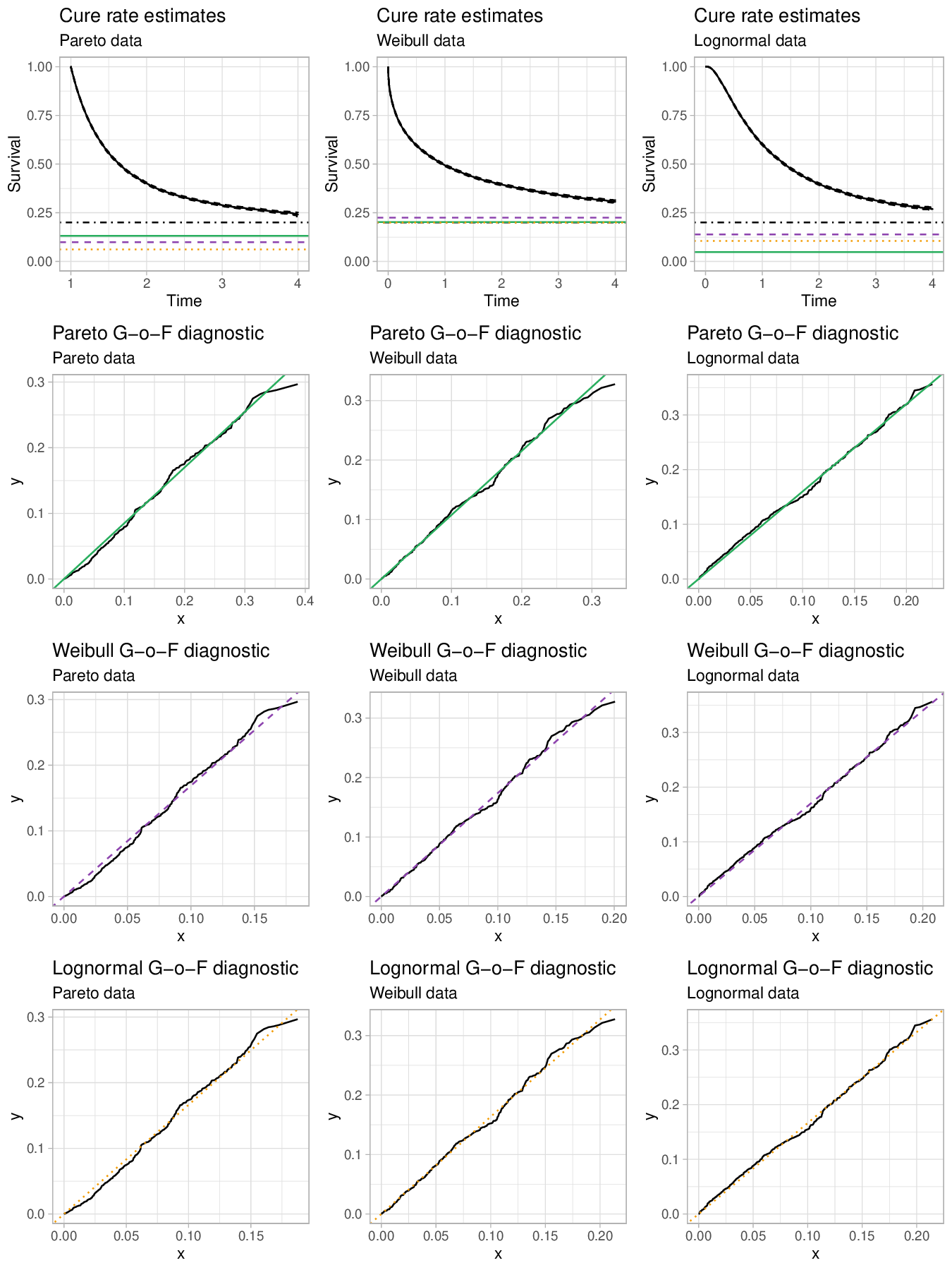}
\caption{Simulated data following the  Pareto, Weibull and lognormal models, best fitting  to the  Norwegian data. Here $k/n=0.1$.}
\label{fig:sims2}
\end{figure}

\begin{figure}[!htbp]
\centering
\includegraphics[width=\textwidth, trim= 0in 0in 0in 0in,clip]{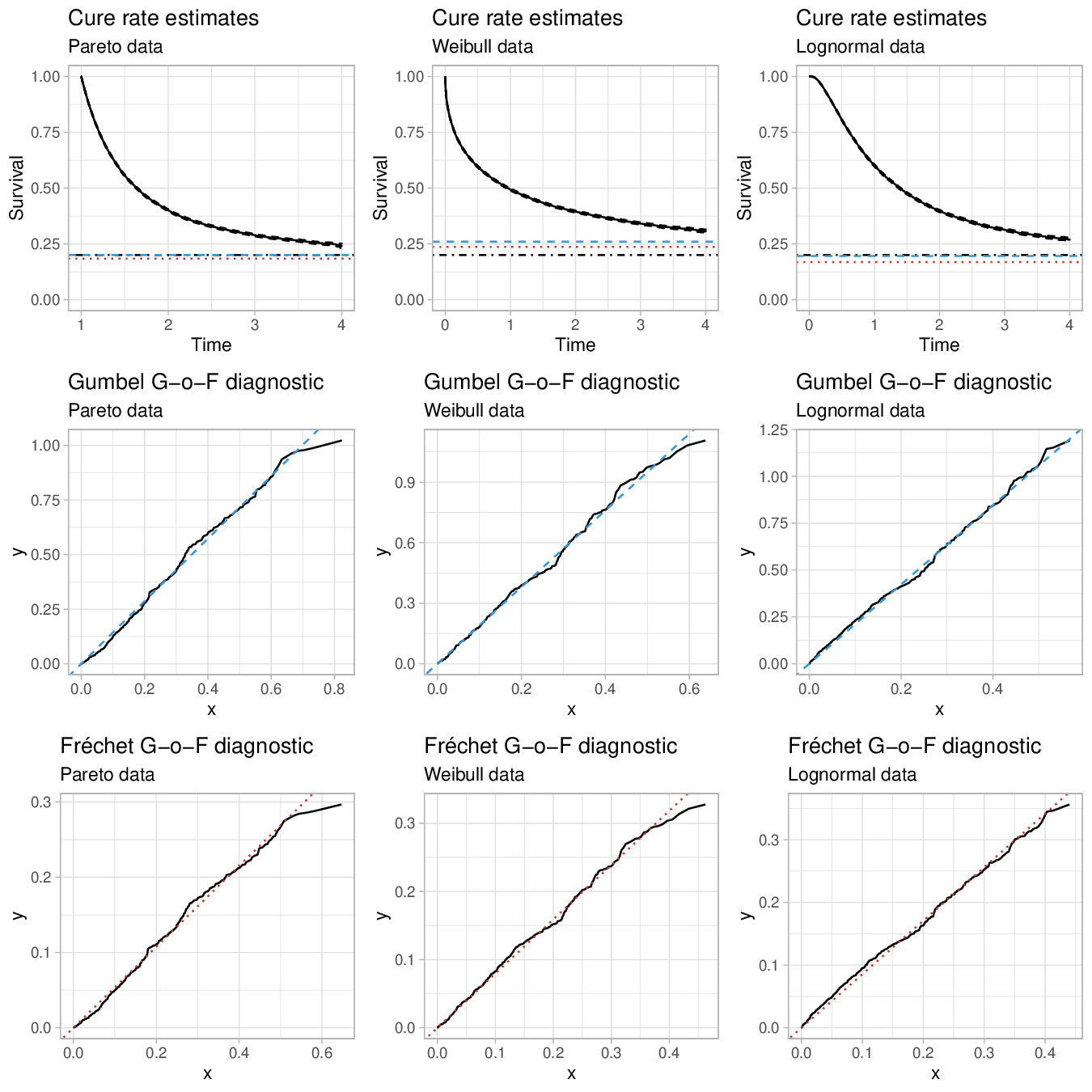}
\caption{Simulated data following the  Pareto, Weibull and lognormal models, best fitting  to the  Norwegian data. Here $k/n=0.1$.}
\label{fig:sims3}
\end{figure}

Next we present the results for $\hat{p}_k^P$, $\hat{p}_k^W$ and $\hat{p}_k^L$ for simulated data sets, comparable to the Norwegian second born data set with respect to sample size ($n=50,000$), using simulated Pareto, Weibull and lognormal distributions for $T$, with parameters as in the simulations study, with $p=0.8$ and censoring distribution uniform on $[0,4]$ in the three cases. In Figure \ref{fig:sims} we have used $k/n=0.9$, while Figure \ref{fig:sims2} corresponds to $k/n=0.1$. We also provide in Figure \ref{fig:sims3} the Gumbel and Fr\'echet approaches, $\hat{p}^G_n$ and $\hat{p}^F_n$, with $k/n=0.1$.
At $k/n=0.9$, the goodness-of-fit plots which correspond to the correct distribution of the data indeed indicate the best fit and lead to satisfactory $p$ estimates. Again for smaller $k/n$ values this link disappears and the different $\hat{p}^*$ estimates become comparable, though naturally more volatile. This appears to be in correspondence with the main conclusions from the asymptotic analysis.   

\begin{figure}[!htbp]
\centering
\includegraphics[width=0.7\textwidth, trim= 0in 0in 0in 0in,clip]{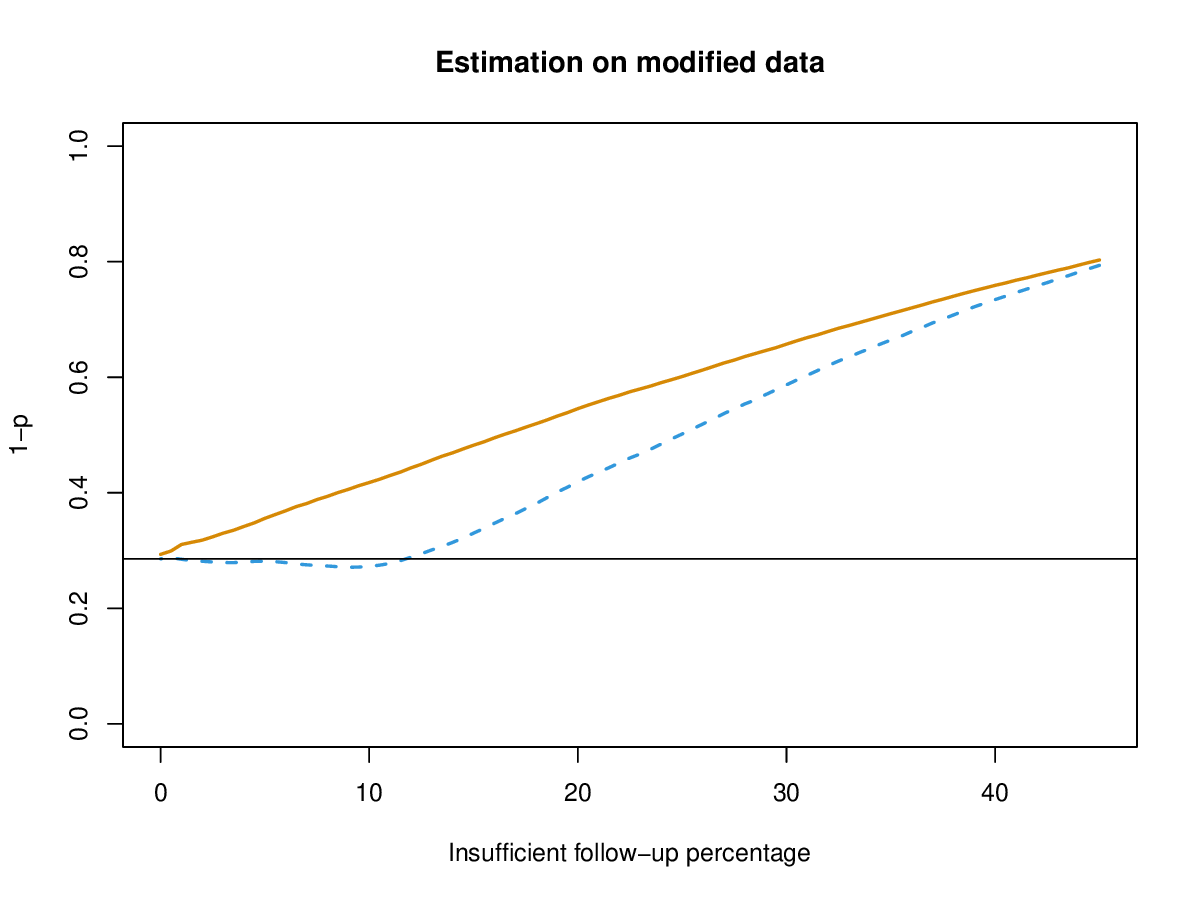}
\caption{Modification of the Norwegian data to increase the insufficient follow-up. We compare the $\hat p_n^G$ (thick dashed) against $p_n$ (thick solid). The $p_G(n,\varepsilon)$ for this specific data is very close to the latter, and not distinguishable in the plot. The horizontal line shows the value of $\hat p_n^G$ for the original (unmodified) data.}
\label{fig:modified_data}
\end{figure}

Finally, consider an alteration of the dataset where we artificially increase the insufficiency of follow-up, by setting observation indicators equal to zero above a certain threshold. This is motivated by the fact that our estimators and the benchmarks provide a very similar estimate, which we know from the simulation study can happen when follow-up is sufficient, and we would like to see a scenario where they depart from each other. The thresholds above which we increase insufficient follow-up are taken as a percentage of the top data, and we consider the range $0-45\%$. We fit $\hat p_n^G$, since it was the most stable estimate across changes of $k$ for the original analysis, and take $k/n=0.5$ (as it should be above $0.45$ to handle the modified data). We observe in Figure \ref{fig:modified_data} that our estimator enjoys stability up to about $12\%$ additional artificial insufficient follow-up, while $p_n$ (and $p_G(n,\varepsilon)$, which for this data follows $p_n$ closely) has no stable region up to any percentage. This suggests that -- when fitting well -- our estimator performs the task it is designed to achieve satisfactorily, and it is rather robust to changes of the censoring distribution upper limit.

\section{Conclusion}\label{sec:conc}
In this paper, we have introduced a non-parametric cure model that integrates extreme value tail estimation to jointly estimate the cure rate and the extreme value index. Our approach uses the full information contained in the top order statistics, improving cure rate estimation in the presence of insufficient follow-up data.

We proposed a Peaks-over-Threshold methodology under the Gumbel max-domain assumption and then extended it to specific models such as Pareto, log-normal, and Weibull tail models. This provides a framework for identifying the most relevant tail characteristics of the susceptible population. Our methods are shown through simulations to rival and often outperform established nonparametric cure rate estimation models in both sufficient and insufficient follow-up scenarios.

Through theoretical asymptotic analysis, we have shown that our estimators maintain desirable weak convergence properties under extreme value conditions, and the regularization mechanism introduced in our probability plotting methodology prevents excessive deviation from standard estimators. While limitations exist, such as sample fraction selection, they present opportunities for further improving the EVT-based estimation techniques for cure models.

\newpage
\appendix

\section{Proof of Theorem~\ref{th:asymptotic}}

The estimating equations based in minimizing \eqref{SSstar} are given by
\begin{align}
&{1 \over k}\sum_{j=1}^k D_{j,k}(\hat{p})\log {Z_{n-j+1,n}\over Z_{n-k,n}} - 
\hat{\beta}_{*,k} {1 \over k}\sum_{j=1}^k\log^2 {Z_{n-j+1,n}\over Z_{n-k,n}} =0 
\label{eq1}
\\
 & {1 \over k}\sum_{j=1}^k 
D_{j,k}(\hat{p})D'_{j,k}(\hat{p}) - \hat{\beta}_{*,k} {1 \over k}\sum_{j=1}^k\left(\log {Z_{n-j+1,n}\over Z_{n-k,n}} \right)D'_{j,k}(\hat{p}) = \lambda (\hat{p}-p_n)
\label{eq2}
\end{align}
with 
\begin{align*}
D_{j,k}(p)&= s_*\left(1- {\hat{F}_n(Z_{n-j+1,n})\over p}\right)- s_*\left(1-{\hat{F}_n(Z_{n-k,n})\over p}\right),\\
D'_{j,k}(p)&= p^{-2}\left(
s'_*\left(1- {\hat{F}_n(Z_{n-j+1,n})\over p}\right)
\hat{F}_n(Z_{n-j+1,n})\right.\\
&\left.-
s'_*\left(1- {\hat{F}_n(Z_{n-k,n})\over p}\right)
\hat{F}_n(Z_{n-k,n})
\right) 
\end{align*}

Using \eqref{decomp} and the mean value theorem we obtain
\begin{align*}
&D_{j,k}(\hat{p}) =\\
&s_*(\bar{F}_{0,*}(Z_{n-j+1,n}))-s_*(\bar{F}_{0,*}(Z_{n-k,n}))\\
&-\left( {p_0(\tau_c) \over \hat{p}}-1 \right) {1 \over F_{0,*}(\tau_c)}
\left[ g_*(\bar{F}_{0,*}(Z_{n-j+1,n}))
 -
 g_*(\bar{F}_{0,*}(Z_{n-k,n}))
 \right](1+o(1))\\
 &  -{1 \over \sqrt{n}p_0(\tau_c)}
 \left[ Y_n (Z_{n-j+1,n}) - Y_n (Z_{n-k,n}) \right] s'_*(\bar{F}_{0,*}(\tau_c))(1+o(1))\\
 & -{\bar{F}_{0,*}(\tau_c)\over F_{0,*}(\tau_c)}
\left[ g_*(\bar{F}_{0,*}(Z_{n-j+1,n}))
 -
 g_*(\bar{F}_{0,*}(Z_{n-k,n}))
 \right](1+o(1)).
 \end{align*}
 With $s^{\leftarrow}_*$ denoting the inverse function of $s_*$, 	 using the mean value theorem we obtain with $g_*(x)=(1-x)s'_*(x)$ that
\begin{align*}
& g_*(\bar{F}_{0,*}(Z_{n-j+1,n}))-g_*(\bar{F}_{0,*}(Z_{n-k,n})) \\
&= \left(
(g_* \circ s^{\leftarrow}_*)(s_*(\bar{F}_{0,*}(Z_{n-j+1,n}))) -
(g_* \circ s^{\leftarrow}_*)(s_*(\bar{F}_{0,*}(Z_{n-k,n})))
 \right) \\
 &= [s_*(\bar{F}_{0,*}(Z_{n-j+1,n})))-(s_*(\bar{F}_{0,*}(Z_{n-k,n}))]A(\tau_c) (1+o(1)).
 \end{align*}  
 From this $D_{j,k}(\hat{p})$ can be further developed as
  \begin{align*}
 &D_{j,k}(\hat{p}) 
 =
s_*(\bar{F}_{0,*}(Z_{n-j+1,n}))-s_*(\bar{F}_{0,*}(Z_{n-k,n}))\\
&-\left( {p_0(\tau_c) \over \hat{p}}-1 \right) {A(\tau_c)(1+o(1)) \over F_{0,*}(\tau_c)}
\left[ s_*(\bar{F}_{0,*}(Z_{n-j+1,n}))
 -
 s_*(\bar{F}_{0,*}(Z_{n-k,n}))
 \right]\\
 &  -{1 \over \sqrt{n}p_0(\tau_c)}
 \left[ Y_n (Z_{n-j+1,n}) - Y_n (Z_{n-k,n}) \right] s'_*(\bar{F}_{0,*}(\tau_c))(1+o(1))\\
 & -{\bar{F}_{0,*}(\tau_c)\over F_{0,*}(\tau_c)}
\left[ s_*(\bar{F}_{0,*}(Z_{n-j+1,n}))
 -
 s_*(\bar{F}_{0,*}(Z_{n-k,n}))
 \right]A(\tau_c)(1+o(1)). 
 \end{align*}

 Similarly we have
\begin{align*}
&D'_{j,k}(\hat{p}) =
s_*(\bar{F}_{0,*}(Z_{n-j+1,n}))-s_*(\bar{F}_{0,*}(Z_{n-k,n}))A(\tau_c)\\
&-\left( {p_0(\tau_c) \over \hat{p}}-1 \right) {B(\tau_c)(1+o(1)) \over F_{0,*}(\tau_c)}
\left[ s_*(\bar{F}_{0,*}(Z_{n-j+1,n}))
 -
 s_*(\bar{F}_{0,*}(Z_{n-k,n}))
 \right]\\
 &  -{1 \over \sqrt{n}p_0(\tau_c)}
 \left[ Y_n (Z_{n-j+1,n}) - Y_n (Z_{n-k,n}) \right] g'_*(\bar{F}_{0,*}(\tau_c))(1+o(1))\\
 & -{\bar{F}_{0,*}(\tau_c)\over F_{0,*}(\tau_c)}
\left[ s_*(\bar{F}_{0,*}(Z_{n-j+1,n}))
 -
 s_*(\bar{F}_{0,*}(Z_{n-k,n}))
 \right]B(\tau_c)(1+o(1)),
  \end{align*}
with $B(\tau_c)= 1-3F_{0,*}(\tau_c)(s''_*/s'_*)(\bar{F}_{0,*}(\tau_c))+  F_{0,*}(\tau_c)(s'''_*/s'_*)(\bar{F}_{0,*}(\tau_c))$.

Based on \eqref{LR} we have for an ordered i.i.d. sequence of $n$ uniform (0,1) order statistics $U_{1,n} \leq U_{2,n}\leq \ldots \leq U_{n,n}$ and similarly $V_{1,k} \leq V_{2,k}\leq \ldots \leq V_{k,k}$ for a sample of size $k$,  as $k,n \to \infty$ and $k/n \to 0$,  
\begin{align*}
 \log Z_{n-j+1,n}-\log Z_{n-k,n} \stackrel{d}{=} \tau_c^{-1}(1-p_0(\tau_c))^{\gamma_c} 
 U_{k+1,n}^{-\gamma_c}\left( 1- V_{j,k}^{-\gamma_c}\right) (1+o_p(1)).
\end{align*}
Similarly
\begin{align*}
&  Z_{n-j+1,n}-Z_{n-k,n}\\
&= (1-p_0(\tau_c))^{\gamma_c}
\left[ -U_{j,n}^{-\gamma_c} (1+O_p(U_{j,n}^{-\gamma_c}))
+ U_{k+1,n}^{-\gamma_c} (1+O_p(U_{k+1,n}^{-\gamma_c}))
\right]
\\
& = (1-p_0(\tau_c))^{\gamma_c} U_{k+1,n}^{-\gamma_c}
\left(1-V_{j,k}^{-\gamma_c}\right)(1+o_p(1)), \\
&  Z_{n-j+1,n}^{-\nu_*}-Z_{n-k,n}^{-\nu_*}\\
& =
-\tau_c^{-\nu_*-1}(1-p_0(\tau_c))^{\gamma_c}\nu_*
U_{k+1,n}^{-\gamma_c}
\left(1-V_{j,k}^{-\gamma_c}\right)(1+o_p(1)).
 \end{align*}
Hence, using \eqref{sF0},
\begin{align*}
&s_*(\bar{F}_{0,*}(Z_{n-j+1,n}))-s_*(\bar{F}_{0,*}(Z_{n-k,n})) \\
&=(\beta_* - D_* \nu_*\tau_c^{-\nu_*})\log {Z_{n-j+1,n} \over Z_{n-k,n}}\\
&= (\beta_* - D_* \nu_*\tau_c^{-\nu_*})\tau_c^{-1}(1-p_0(\tau_c))^{\gamma_c} U_{k+1,n}^{-\gamma_c}
\left(1-V_{j,k}^{-\gamma_c} \right)(1+o_p(1)).
\end{align*}

The first equation \eqref{eq1} can now be rewritten as
\begin{align*}
& \left([\beta_* -D_*\nu_*\tau_c^{-\nu_*}]- \hat{\beta}_*\right) 
{1 \over k}\sum_{j=1}^k \log^2 {Z_{n-j+1,n} \over Z_{n-k,n}}\\
& - A(\tau_c)\left( { p_0(\tau_c)\over \hat{p}_k } -1\right)
{1 \over F_{0,*}(\tau_c)}[\beta_* -D_*\nu_*\tau_c^{-\nu_*}]
 {1 \over k}\sum_{j=1}^k \log^2 {Z_{n-j+1,n} \over Z_{n-k,n}}\\
 &
 =
 {s'_*(\bar{F}_{0,*}(\tau_c))\over  p_0(\tau_c)(1-p_0(\tau_c))} \tilde{T} _{k,n}
\\
 &  +A(\tau_c){\bar{F}_{0,*}(\tau_c) \over F_{0,*}(\tau_c)} [\beta_* -D_*\nu_*\tau_c^{-\nu_*}]
 {1 \over k}\sum_{j=1}^k \log^2 {Z_{n-j+1,n} \over Z_{n-k,n}}, 
  \end{align*}
 with $$\tilde{T} _{k,n}=  {1 \over  \sqrt{n}k}\sum_{j=1}^k \log {Z_{n-j+1,n} \over Z_{n-k,n}}\left[ Y_n (Z_{n-j+1,n}) - Y_n (Z_{n-k,n}) \right]. $$
Next we use the above approximation of $\log (Z_{n-j+1,n}/Z_{n-k,n})$ and that $${1 \over k}\sum_{j=1}^k \left(1-V_{j,k}^{-\gamma_c}\right)^2=  {2\gamma_c^2 \over (1-\gamma_c)(1-2\gamma_c)}(1+o_p(1))$$ as $k \to \infty$ and  $kU_{k+1,n}/n = 1+o_p(1))$ as $k,n \to \infty$. 

Moreover we approximate $Y_n(Z_{n-j+1,n})$ by $(1-F)(Z_{n-j+1,n})\mathbf{Z}(Z_{n-j+1,n})$  ($j=1,...,k$) based on the almost sure convergence of $Y_n$ to $(1-F)\mathbf{Z}$ uniformly on $[0,\tau_c]$. 
 
Next we approximate  $\mathbf{Z}(Z_{n-j+1,n})$ by $\mathbf{Z}(U_H({n+1 \over j}))$, using Theorem 2.4.2 in~\cite{haan2006extreme} based on \eqref{LR}, obtaining
 \begin{align}
 Z_{n-j+1,n} &= U_H ({n+1 \over j}) 
 + {1 \over\sqrt{k}} a({n \over k}) W({j \over k+1}) ({j \over k+1})^{-\gamma_c-1}(1+o_p(1)),
 \label{ZU}
 \end{align}
for $j=1, \ldots,k,$  where $a(n/k)= -\gamma_c(n/k)^{\gamma_c}(1-p_0(\tau_c))$, $W$ a Brownian motion, and $o_p(1)$ holding uniformly in $j=1,\ldots,k$.
 Furthermore note that
 the variance  function $v$  with $t \to \infty$ and $x$ bounded satisfies
\begin{align}
& v(U_H(t+x))- v(U_H(t)) \nonumber \\&= p\int_{U_H(t)}^{U_H(t+x)}{ \dd F_0(s) \over \bar{F}(s)(1-H(s))} \nonumber\\
& = p \int_t^{t+x} {u \over \bar{F}(U_H(u))}\dd F_{0,*}(U_H(u)) \nonumber \\
&= p  \int_t^{t+x} {u \over 1-p_0 (\tau_c) (1+o(1))}
\dd F_{0,*}(\tau_c - (1-p_0(\tau_c))^{\gamma_c}u^{\gamma_c})(1+o(1)) \nonumber \\
&= p  \int_t^{t+x} {u \over 1-p_0 (\tau_c) (1+o(1))}
\dd\left\{F_{0,*}(\tau_c)-(1-p_0(\tau_c))^{\gamma_c}u^{\gamma_c}f_{0,*}(\tau_c) (1+o(1))\right\} \nonumber \\
&= 
 B_v x t^{\gamma_c}(1+o(1)), \label{dv1}
\end{align}
 with $B_v= p\gamma_c (1-p_0(\tau_c))^{\gamma_c-1} f_{0,*}(\tau_c)$. 
 Writing ${\mathbf Z}(t)={\mathbf B} (v(t))$ with $ {\mathbf B}$ a Brownian motion, combining \eqref{ZU} and \eqref{dv1} we obtain
 \begin{align*}
 & \mathbf{Z}(Z_{n-j+1,n}) \\&=  {\mathbf B} \left( v \left[ 
  U_H ({n+1 \over j})+ 
  {(1+o_p(1)) \over\sqrt{k}} a({n \over k}) W({j \over k+1}) ({j \over k+1})^{-\gamma_c-1}
 \right] \right) \\
 &= {\mathbf B} \left( 
 v( U_H ({n+1 \over j}))\right)\\
 &\quad+  B_v {1 \over\sqrt{k}} a({n \over k}) \left(U_H({n+1 \over j})\right)^{\gamma_c}
   W({j \over k+1}) ({j \over k+1})^{-\gamma_c-1}(1+o_p(1))  \\
   &= {\mathbf Z} ( U_H ({n+1 \over j})) + O_p(k^{-1/4}
a^{1/2}(n/k)), \; \mbox{ uniformly in } j=1,\ldots,k, \end{align*}
where the last line follows from Lemma 2.2 in~\cite{BEIRLANT1990241}.

Now the first equation is asymptotically equivalent to
\begin{align*}
& \left([\beta_* -D_*\nu_*\tau_c^{-\nu_*}]- \hat{\beta}_*\right) \left({k \over n}\right)^{-2\gamma_c}{1 \over \tau^2_c}
(1-p_0(\tau_c))^{2\gamma_c}{2\gamma_c^2 \over (1-\gamma_c)(1-2\gamma_c)}
\\
& - \left( { p_0(\tau_c)\over \hat{p}_k } -1\right)\left({k \over n}\right)^{-2\gamma_c}
{A(\tau_c) \over F_{0,*}(\tau_c)}[\beta_* -D_*\nu_*\tau_c^{-\nu_*}]{1 \over \tau^2_c}
{2\gamma_c^2 (1-p_0(\tau_c))^{2\gamma_c} \over (1-\gamma_c)(1-2\gamma_c)}
\\
&=\left({k \over n}\right)^{-\gamma_c} T_{k,n}  {s'_*(\bar{F}_{0,*}(\tau_c))\over p_0(\tau_c)} {1 \over \tau_c}
(1-p_0(\tau_c))^{\gamma_c} \\
&\quad +\left({k \over n}\right)^{-2\gamma_c}A(\tau_c){\bar{F}_{0,*}(\tau_c) \over F_{0,*}(\tau_c)}
[\beta_* -D_*\nu_*\tau_c^{-\nu_*}]
{1 \over \tau^2_c}
(1-p_0(\tau_c))^{2\gamma_c}{2\gamma_c^2 \over (1-\gamma_c)(1-2\gamma_c)},
\end{align*}
which directly leads to the final version of the first equation. 

In order to simplify the second equation \eqref{eq2}, note that 
\begin{align*}
& \hat{p} - \hat{F}_n (Z_{n,n}) \\
 &= 
(\hat{p} - p_0(\tau_c)) -n^{-1/2}Y_n (Z_{n,n}) 
- p \left( F_{0,*}(U_H(U_{1,n}^{-1}))-F_{0,*}(\tau_c)   \right)\\
&=(\hat{p} - p_0(\tau_c)) -n^{-1/2}Y_n (Z_{n,n})  
+ U_{1,n}^{-\gamma_c}(1-p_0(\tau_c))^{\gamma_c}f_{0,*}(\tau_c)(1+o(1)). 
\end{align*}
Now we obtain similarly for  the second equation \eqref{eq2} 
\begin{align*}
&{1 \over k}\sum_{j=1}^k \Big\{ [s_*(\bar{F}_{0,*}(Z_{n-j+1,n}))-s_*(\bar{F}_{0,*}(Z_{n-k,n}))]
\\
& 
+ {A(\tau_c) \over F_{0,*}(\tau_c)}[s_*(\bar{F}_{0,*}(Z_{n-j+1,n}))-s_*(\bar{F}_{0,*}(Z_{n-k,n}))]
\left(1-{p_0(\tau_c)\over \hat{p}_{*,k}}\right)\\
& 
- {s'_*(\bar{F}_{0,*}(\tau_c))\over \sqrt{n}p_0(\tau_c) }[Y_n(Z_{n-j+1,n})-Y_n(Z_{n-k,n})]\\
& 
- A(\tau_c){\bar{F}_{0,*}(\tau_c) \over F_{0,*}(\tau_c)}[s_*(\bar{F}_{0,*}(Z_{n-j+1,n}))-s_*(\bar{F}_{0,*}(Z_{n-k,n}))]\Big\}\\
&
\times \Big\{ [s_*(\bar{F}_{0,*}(Z_{n-j+1,n}))-s_*(\bar{F}_{0,*}(Z_{n-k,n}))] A(\tau_c)
\\
&
+ {B(\tau_c) \over F_{0,*}(\tau_c)}[s_*(\bar{F}_{0,*}(Z_{n-j+1,n}))-s_*(\bar{F}_{0,*}(Z_{n-k,n}))]
\left(1-{p_0(\tau_c)\over \hat{p}_{*,k}}\right)\\
& 
- {g'_*(\bar{F}_{0,*}(\tau_c))\over \sqrt{n}p_0(\tau_c) }[Y_n(Z_{n-j+1,n})-Y_n(Z_{n-k,n})]\\
&
- B(\tau_c){\bar{F}_{0,*}(\tau_c) \over F_{0,*}(\tau_c)}[s_*(\bar{F}_{0,*}(Z_{n-j+1,n}))-s_*(\bar{F}_{0,*}(Z_{n-k,n}))]\Big\}\\
& - \hat{\beta}_{*,k}
{1 \over k}\sum_{j=1}^k \log {Z_{n-j+1,n} \over Z_{n-k,n}} \Big\{
[s_*(\bar{F}_{0,*}(Z_{n-j+1,n}))-s_*(\bar{F}_{0,*}(Z_{n-k,n}))] A(\tau_c) \\ 
&
 + {B(\tau_c) \over F_{0,*}(\tau_c)}[s_*(\bar{F}_{0,*}(Z_{n-j+1,n}))-s_*(\bar{F}_{0,*}(Z_{n-k,n}))]
\left(1-{p_0(\tau_c)\over \hat{p}_{*,k}}\right)\\
& 
- {g'_*(\bar{F}_{0,*}(\tau_c))\over \sqrt{n}p_0(\tau_c) }[Y_n(Z_{n-j+1,n})-Y_n(Z_{n-k,n})] \\
& -  B(\tau_c){\bar{F}_{0,*}(\tau_c) \over F_{0,*}(\tau_c)}[s_*(\bar{F}_{0,*}(Z_{n-j+1,n}))-s_*(\bar{F}_{0,*}(Z_{n-k,n}))]\Big\} \\
&
= {\lambda \over p}\left( \hat{p}_{*,k}- p_0(\tau_c)\right) -{\lambda \over \sqrt{n}p}Y_n (Z_{n,n})+
\lambda U_{1,n}^{-\gamma_c}(1-p_0(\tau_c))^{\gamma_c}
f_{0,*}(\tau_c),
\end{align*}
or, when replacing  $\hat{\beta}_{*,k}$ by $\beta_* -D_*\nu_* \tau_c^{-\nu_*}$ and coupled with the last three terms in $D'_{j,k}(\tau_c)$,
\begin{align*}
& \left([\beta_* -D_*\nu_*\tau_c^{-\nu_*}]- \hat{\beta}_*\right) 
{1 \over k}\sum_{j=1}^k \log^2 {Z_{n-j+1,n} \over Z_{n-k,n}}A(\tau_c)[\beta_* -D_*\nu_*\tau_c^{-\nu_*}]\\
&+\left(\hat{p}_{*,k}-p_0 (\tau_c)\right){1 \over p_0(\tau_c)}
\\
&
\times\left\{
-\lambda {p_0(\tau_c) \over p} +[\beta_* -D_*\nu_*\tau_c^{-\nu_*}]^2
{1 \over k}\sum_{j=1}^k \log^2 {Z_{n-j+1,n} \over Z_{n-k,n}}{A^2(\tau_c) \over F_{0,*}(\tau_c)}
\right\} \\
&= 
{A(\tau_c) \over p_0(\tau_c) }[\beta_* -D_*\nu_*\tau_c^{-\nu_*}]s'_*(\bar{F}_{0,*}(\tau_c))
\;\\
&\times{1 \over \sqrt{n}k}\sum_{j=1}^k \log {Z_{n-j+1,n} \over Z_{n-k,n}}[Y_n(Z_{n-j+1,n})-Y_n(Z_{n-k,n})]
\\
&  -{\lambda \over p} {1 \over \sqrt{n}}Y_n(Z_{n,n}) \\
& +A^2(\tau_c){\bar{F}_{0,*}(\tau_c) \over F_{0,*}(\tau_c)}[\beta_* -D_*\nu_*\tau_c^{-\nu_*}]^2\\
&\times
{1 \over k}\sum_{j=1}^k \log^2 {Z_{n-j+1,n} \over Z_{n-k,n}} +\lambda U_{1,n}^{-\gamma_c}(1-p_0(\tau_c))^{\gamma_c} f_{0,*}(\tau_c).
\end{align*}
With the use of Lemma 1 in~\cite{EscobarKeilegom2019} this equation is asymptotically equivalent to
\begin{align*}
& \left([\beta_* -D_*\nu_*\tau_c^{-\nu_*}]- \hat{\beta}_*\right) 
A(\tau_c)[\beta_* -D_*\nu_*\tau_c^{-\nu_*}]
\left({k \over n}\right)^{-2\gamma_c}{(1-p_0(\tau_c))^{\gamma_c}\over \tau_c}
{1 \over M(\tau_c)}
\\
& +\frac{\hat{p}_{*,k}-p_0 (\tau_c)}{ p_0(\tau_c)}\\
&\times
\left\{
-\lambda {p_0(\tau_c) \over p} +[\beta_* -D_*\nu_*\tau_c^{-\nu_*}]^2
{A^2(\tau_c) \over F_{0,*}(\tau_c)}
 \left({k \over n}\right)^{-2\gamma_c}{(1-p_0(\tau_c))^{\gamma_c}\over \tau_c}
{1 \over M(\tau_c)}\right\}\\
&= \left({k \over n}\right)^{-\gamma_c}T_{k,n}
{A(\tau_c) \over p_0(\tau_c) }[\beta_* -D_*\nu_*\tau_c^{-\nu_*}]s'_*(\bar{F}_{0,*}(\tau_c))
{(1-p_0(\tau_c))^{\gamma_c} \over \tau_c}
\\
&  -{\lambda \over p} {1 \over \sqrt{n}}Y_n(U_H(n)) \\
&  
+A^2(\tau_c){\bar{F}_{0,*}(\tau_c) \over F_{0,*}(\tau_c)}[\beta_* -D_*\nu_*\tau_c^{-\nu_*}]^2
\left({k \over n}\right)^{-2\gamma_c}{(1-p_0(\tau_c))^{\gamma_c}\over \tau_c}
{1 \over M(\tau_c)}
 \\
& +\lambda U_{1,n}^{-\gamma_c}(1-p_0(\tau_c))^{\gamma_c} f_{0,*}(\tau_c),
\end{align*}
which leads to the stated version of the second equation when dividing by 
$$[\beta_* -D_*\nu_*\tau_c^{-\nu_*}]
\left({k \over n}\right)^{-2\gamma_c}{(1-p_0(\tau_c))^{\gamma_c}\over \tau_c}
{1 \over M(\tau_c)}$$ and observing that $$\lambda_{k,n} U_{1,n}^{-\gamma_c}= O_p ((k/n)^{-2\gamma_c}n^{\gamma_c})$$ is asymptotically negligible compared to $\lambda_{k,n}n^{-1/2}Y_n(Z_{n,n})$.

Similarly as in \eqref{dv1}, we obtain
\begin{align}
 v(U_H(tx))- v(U_H(t)) = B_v x^{1+\gamma_c} h_{1+\gamma_c}(t)(1+o(1)).
\label{dv2}
\end{align}
For the variance of $(n/k)^{-\gamma_c} T_{k,n}$
 we then find
 \begin{align*}
&  (1-p_0(\tau_c))^2\left({n \over k}\right)^{-2\gamma_c}n^{-1}
 {1 \over k^2}
\sum_{j_1=1}^k \sum_{j_2=1}^k 
\left(1- \left( {j_1 \over k+1}\right)^{-\gamma_c} \right) 
\left(1- \left( {j_2 \over k+1}\right)^{-\gamma_c} \right)\\
& \times
\E\left[
[\mathbf{Z}(U_H({n+1 \over j_1}))- \mathbf{Z}(U_H({n+1 \over k+1}))][\mathbf{Z}(U_H({n+1 \over j_2}))- \mathbf{Z}(U_H({n+1 \over k+1}))]
\right] \\
&= B_v (1-p_0(\tau_c))^2 \left({n \over k}\right)^{-2\gamma_c}n^{-1}
[h_{1+\gamma_c}\left( {n+1 \over j_1 \vee j_2} \right)- h_{1+\gamma_c}\left( {n+1 \over k+1} \right)]\\
& \times
{1 \over k^2}
\sum_{j_1=1}^k \sum_{j_2=1}^k 
\left(1- \left( {j_1 \over k+1}\right)^{-\gamma_c} \right) 
\left(1- \left( {j_2 \over k+1}\right)^{-\gamma_c} \right)
 \\
&= B_v (1-p_0(\tau_c))^2 \left({n \over k}\right)^{-2\gamma_c}n^{-1}
\left({n \over k}\right)^{\gamma_c+1}\sigma^2_k
= B_v (1-p_0(\tau_c))^2  \left({n \over k}\right)^{-\gamma_c}k^{-1}\sigma^2_k.
 \end{align*}
Similarly, for the correlation between the terms $(n/k)^{-\gamma_c}(1-p_0(\tau_c))^{-1}T_{k,n}$ and $n^{-1/2}\mathbf{Z}(U_H(n+1))$ we obtain
  \begin{align*}
 & \left( {n \over k}\right)^{-\gamma_c}n^{-1} 
  {1 \over k}\sum_{j=1}^k \left( 1- \left( {j \over k+1}\right)^{-\gamma_c}\right)
  \\
  &\times\E\left[ \left[\mathbf{Z}(U_H({n+1 \over j}))- \mathbf{Z}(U_H({n+1 \over k+1})) \right]\mathbf{Z}(U_H(n+1))\right]\\
  &= 
  B_v \left( {n \over k}\right)^{-\gamma_c}{1 \over nk}\sum_{j=1}^k \left( 1- \left( {j \over k+1}\right)^{-\gamma_c}\right)\left[ h_{1+\gamma_c}\left( {n+1 \over j}\right)- h_{1+\gamma_c}\left( {n+1 \over k+1}\right)\right]\\
 &= B_v \left( {n \over k}\right)^{-\gamma_c}n^{-1} \left( {n \over k}\right)^{1+\gamma_c}{1 \over k}\sum_{j=1}^k \left( 1- \left( {j \over k+1}\right)^{-\gamma_c}\right)h_{1+\gamma_c}\left( {k+1 \over j} \right)\\ &= O\left( {1 \over k} \right)
 = o\left( k^{-1} \left({n \over k} \right)^{-\gamma_c}\right).
  \end{align*}

This finishes the proof.

\newpage
\bibliographystyle{imsart-nameyear} 
\bibliography{ejs-sample.bib}       

\end{document}